\input ./style/arxiv-general.cfg
\documentclass[aos,MSNbibl,dvips]{arximspdf}
\makeatletter
   \@ifpackageloaded{graphicx}{}{\usepackage{graphicx}}
\makeatother
\usepackage{stfloats}

\doi{10.1214/14-AOS1300}
\volume{43}
\issue{3}
\pubyear{2015}
\firstpage{1089}
\lastpage{1116}
\docsubty{FLA}

\makeatletter
  \fnbelowfloat
\def\sfrac#1#2{#1/#2}

\def\sklfrac#1#2{(#1/#2)}

\newcommand{\rrvert}{\vert}
\newcommand{\rrVert}{\Vert}
\newcommand{\llvert}{\vert}
\newcommand{\llVert}{\Vert}
\renewcommand{\mid}{|}
\newtheorem{theorem}{Theorem}
\newtheorem{hypo}{Hypothesis}
\newtheorem{lemma}{Lemma}
\newproclaim{definition}{Definition}
\newproclaim{remark}{Remark}
\newcommand{\binom}[2]{{#1 \choose #2}}
\newcommand{\reals}{{\mathbb{R}}}
\newcommand{\integers}{{\mathbb{Z}}}
\newcommand{\naturals}{{\mathbb{N}}}
\newcommand{\eexp}{\mathrm{e}}
\newcommand{\diff}{\mathrm{d}}
\newcommand{\Expect}{\mathbb{E}}
\newcommand{\Prob}{\mathbb{P}}
\newcommand{\Indc}{\mathbf{1}}
\newcommand{\bbP}{{\mathbb{P}}}
\newcommand{\sfX}{{\mathsf{X}}}
\newcommand{\sfY}{{\mathsf{Y}}}
\newcommand{\calE}{{\mathcal{E}}}
\newcommand{\calF}{{\mathcal{F}}}
\newcommand{\calG}{{\mathcal{G}}}
\newcommand{\calL}{{\mathcal{L}}}
\newcommand{\calM}{{\mathcal{M}}}
\newcommand{\calP}{{\mathcal{P}}}
\newcommand{\calQ}{{\mathcal{Q}}}
\newcommand{\bfA}{{\mathbf{A}}}
\newcommand{\bfX}{{\mathbf{X}}}
\newcommand{\tB}{{\widetilde{B}}}
\newcommand{\tX}{{\widetilde{X}}}
\newcommand{\TV}{\mathsf{TV}}
\newcommand{\RR}{\mathbb{R}}
\newcommand{\tr}{\mathsf{Tr}}
\renewcommand{\bfA}{A}
\renewcommand{\bfX}{X}
\newcommand{\bX}{\breve{X}}
\newcommand{\PC}{\mathsf{PC}}
\newcommand{\Tlin}{T_{\mathrm{lin}}}
\newcommand{\Tscan}{T_{\mathrm{scan}}}
\newcommand{\Tmax}{T_{\max}}
\newcommand{\quana}[3][]{[ #3 ]_{#2}^{#1}}
\newcommand{\tquana}[3][]{[#3]_{#2}^{#1}}

\makeatother

\begin{document}
\begin{frontmatter}

\title{Computational barriers in minimax submatrix~detection}
\runtitle{Submatrix detection}

\begin{aug}
\author[A]{\fnms{Zongming}~\snm{Ma}\corref{}\thanksref{T1}\ead
[label=e2]{zongming@wharton.upenn.edu}}
\and
\author[B]{\fnms{Yihong}~\snm{Wu}\ead[label=e3]{yihongwu@illinois.edu}}
\runauthor{Z.~Ma and Y.~Wu}
\affiliation{University of Pennsylvania and University of Illinois at
Urbana-Champaign}
\address[A]{Department of Statistics\\
The Wharton School\\
University of Pennsylvania\\
Philadelphia, Pennsylvania 19104\\
USA\\
\printead{e2}}
\address[B]{Department of Electrical\\
\quad and Computer Engineering\\
University of Illinois at Urbana-Champaign\\
Urbana, Illinois 61801\\
USA\\
\printead{e3}}
\end{aug}
\thankstext{T1}{Supported in part by NSF Career Award DMS-13-52060 and
the Claude Marion Endowed Faculty Scholar Award of the Wharton School.}

%
\received{\smonth{8} \syear{2014}}
%
\revised{\smonth{12} \syear{2014}}

%
\begin{abstract}
This paper studies the minimax detection of a small submatrix of
elevated mean
in a large matrix contaminated by additive Gaussian noise.
To investigate the tradeoff between statistical performance and computational
cost from a complexity-theoretic perspective, we consider a sequence of
discretized models which are asymptotically equivalent to the Gaussian model.
Under the hypothesis that the planted clique detection problem cannot be
solved in randomized polynomial time when the clique size is of smaller order
than the square root of the graph size, the following phase transition
phenomenon
is established: when the size of the large matrix $p \to\infty$, if the
submatrix size $k = \Theta(p^\alpha)$ for any $\alpha\in(0, {2}/{3})$,
computational complexity constraints can incur a severe penalty on the
statistical performance in the sense that any randomized
polynomial-time test is minimax suboptimal by a polynomial factor in $p$;
if $k = \Theta(p^\alpha)$ for any $\alpha\in({2}/{3}, 1)$, minimax optimal
detection can be attained within constant factors in linear time.
Using Schatten norm loss as a representative example, we show that the
hardness of attaining the minimax estimation rate can crucially depend
on the loss function.
Implications on the hardness of support recovery are also obtained.
\end{abstract}

%
\end{frontmatter}

\setcounter{footnote}{1}
\section{Introduction}\label{secintro}

Statistical inference of structured large matrices lies at the heart of
many applications involving massive datasets, such as matrix
completion, functional genomics, community detection and clustering;
see, for instance, \mbox{\cite{Candes09,Shabalin09,BI12,Arias13a,Verzelen13}}
and the references therein.
Many of these detection and estimation problems have been investigated
from a decision theoretic viewpoint,
where
one first establishes a minimax lower bound for any test or estimator
and then constructs a specific
procedure which attains 
the lower bound within a constant or logarithmic factor.

An important element absent from the foregoing decision theoretic\break 
paradigm is computational complexity.
This aspect is especially relevant in the context of high-dimensional
statistical inference, where computationally efficient procedures
(e.g., convex programming, iterative algorithms, etc.) are highly desirable.
However, it has been empirically observed in several basic detection
and estimation problems that popular low-complexity algorithms fail to
attain the minimax rates; see, for example, \cite
{BKRSW11,BI12,Berthet12,Berthet13,Arias13a,KNV13}.
This invites the following question:
how much do we need to back off from the statistical optimality due to
computational complexity constraints?
In this paper, we revisit the sparse submatrix detection problem
that has been studied in \cite
{Shabalin09,KBRS11,BKRSW11,BI12,Sun13,Bhamidi12}, where the goal is to
detect a small submatrix with elevated mean in a large noisy matrix.
Motivations for this detection problem include biclustering for
analyzing microarray data \cite{Shabalin09} and community detection in
social networks \cite{Arias13a}, etc.

\subsection{Problem formulation}
Let\vspace*{-1pt} $X = (X_{ij})$ be a $p\times p$ matrix with independent Gaussian
entries $X_{ij} \stackrel{\mathrm{ind.}}{\sim} N(\theta_{ij},1)$.
Denote the mean matrix by $\theta= (\theta_{ij} ) \in\reals^{p
\times p}$ and the distribution of $X$ by $\Prob_\theta$.
The submatrix detection problem deals with the following setup \cite{BI12}:
under the null hypothesis, the signal is absent, and $\theta$ is a
zero matrix.
Under the alternative hypothesis, $\theta$ is zero except for a
submatrix of size at least $k\times k$ where all the entries exceed
some positive value $\lambda$.
In other words, detecting the submatrix boils down to testing the
following hypotheses on the mean matrix:
%
\begin{equation}
H_0\dvtx X \sim\Prob_0\quad\mbox{versus}\quad
H_1\dvtx X \sim\Prob_\theta,\qquad \theta\in\calM(p,k,\lambda),
\label{eqHgauss}
\end{equation}
where $\Prob_0$ is standard Gaussian, and the parameter space for the
alternative hypothesis is
%
\begin{eqnarray}\label{eqset-H1+}
\calM(p,k,\lambda) &= & \bigl\{\theta\in
\reals^{p\times p}\dvtx \exists U,V\subset[p],\mbox{ s.t.
}\llvert U\rrvert,
\llvert V\rrvert\geq k,
\nonumber\\[-8pt]\\[-8pt]\nonumber
&&\hspace*{5pt} \theta_{ij} \geq\lambda,\mbox{ if } (i,j)\in U\times V,
\theta_{ij} = 0 \mbox{ if }(i,j)\notin U\times V \bigr\}.
\end{eqnarray}
In this problem, the key parameters are the matrix dimension $p$, the
block size $k$ and the signal magnitude $\lambda$. Clearly, it is
easier to detect the submatrix if either $k$ or $\lambda$ increases.
Throughout the paper, we focus on the asymptotic setting where $p$
tends to infinity and both $k = k(p)$ and $\lambda= \lambda(p)$ are
functions of $p$, though we typically drop the explicit dependence on
$p$ for conciseness.

For any test $\phi\dvtx \reals^{p\times p} \to\{0,1\}$, we denote its
worst-case Type-\textup{I${}+{}$II} error probability of testing (\ref{eqHgauss}) by
%
\begin{equation}
\label{eqperr-phi} \calE(\phi) = \Prob_0 \bigl\{ \phi(X) = 1 \bigr\} +
\sup_{\theta
\in\calM
(p,k,\lambda)}\Prob_\theta\bigl\{ \phi(X) = 0 \bigr\}.
\end{equation}
The optimal total probability of error is denoted by
%
\begin{equation}
\calE^* = \inf_{\phi\dvtx \reals^{p\times p}\to\{ 0,1 \}} \calE(\phi).
\label{eqperr}
\end{equation}
In the asymptotic regime of
%
\begin{eqnarray}
\label{eqasymp}
&& p\to\infty,\qquad k\to\infty\quad\mbox{and}\quad k/p\to0,
\end{eqnarray}
the necessary and sufficient condition for reliably detecting the
submatrix has been characterized by Butucea and Ingster (\cite{BI12},
Theorems~2.1 and 2.2):
$\calE^* \to0$ if
%
\begin{equation}
\frac{\lambda}{p/k^2}\to\infty\quad\mbox{or}\quad\liminf_{p\to
\infty
}
\frac{\lambda}{2\sqrt{({1}/{k})\log ({p}/{k})}} > 1, \label
{eqminimax-possible}
\end{equation}
and, conversely, $\calE^* \to1$ if
%
\begin{equation}
\frac{\lambda}{p/k^2}\to0\quad\mbox{and}\quad\limsup_{p\to\infty}
\frac{\lambda}{2\sqrt{({1}/{k})\log({p}/{k})}} < 1. \label
{eqminimax-impossible}
\end{equation}
From this point forward, we say reliable detection is statistically
possible if $\calE^*\to0$ and a sequence of tests $\{\phi_p\}$
reliably detects the submatrix if $\calE(\phi_p)\to0$.

To reliably detect the submatrix under condition (\ref
{eqminimax-possible}), Butucea and Ingster \cite{BI12} proposed a test
involving enumerating all $k\times k$ submatrices of $X$, which is
asymptotically optimal but computationally intensive. It is unclear
from first principles whether statistically optimal detection can be
achieved using computationally efficient procedures.
Thus an intriguing question is in order: under the optimal condition
(\ref{eqminimax-possible}) so that $\calE^*\to0$,
is there a sequence of computationally efficient tests $\{\phi_p\}$
such that $\calE(\phi_p)\to0$?

\subsection{The penalty incurred by complexity constraints}

To approach the computational hardness of the submatrix detection
problem rigorously, we need to investigate the computational cost of
testing procedures in a complexity theoretic sense. However, an
immediate hurdle for the Gaussian model (\ref{eqHgauss}) is that
computational complexity is not well defined for all tests dealing with
nondiscrete distributions
since the observation cannot be represented by finitely many bits
almost surely.
To propose a paradigm for complexity-constrained hypotheses testing, we
consider a sequence of \emph{discretized} Gaussian models which is
asymptotically equivalent to the original model in the sense of Le Cam
\cite{LeCam86} and hence preserves the statistical difficulty of the problem.
More importantly, the computational complexity of tests on the
discretized model can be appropriately defined. See Section~\ref
{secdiscrete} for details.

Next, we take the standard reduction approach in complexity theory: we
show that if the signal magnitude is smaller
than a certain threshold, detecting the submatrix is computationally no
easier than certain well-known intractable problems.
In other words, if an efficient method existed for submatrix detection,
it would lead to an efficient solution to this problem.
In this paper, we use the \emph{planted clique} problem as the
benchmark, which deals with detecting whether a given instance of an
Erd\H{o}s--R\'enyi random graph of size $N$ contains a planted clique
of size $\kappa$.
It is widely believed that the detection problem cannot be solved in
randomized polynomial time when $\kappa= o(\sqrt{N})$,
which we shall refer to as the \emph{planted clique hypothesis}.
For the precise statement and further discussions, see Definition~\ref
{defpc} and Hypothesis~\ref{hypopc} in Section~\ref{seccomplexity}.

Assuming the planted clique hypothesis, our main finding (Theorem~\ref{thmad} in Section~\ref{seccomplexity}) characterizes when it is
possible to achieve reliable detection using computationally efficient
procedures and when it is impossible.
The core of the arguments lies in a \emph{randomized polynomial-time
reduction} scheme which maps
the $N\times N$ adjacency matrix of the random graph in the planted
clique problem to
a $p\times p$ random matrix in polynomial time.
It is worth noting that when $k\geq p^\alpha$ for some $\alpha\geq
\frac{1}{2}$, the cardinality of the graph $N$
is not equal to the size of the matrix $p$ but rather chosen to be
$p^{1+\delta}$ (omitting log factor), where $\delta> 0$ depends on
$\alpha$.
On the other hand, $\kappa$ can always be chosen as a constant
multiple of $k$.

\begin{figure}

\includegraphics{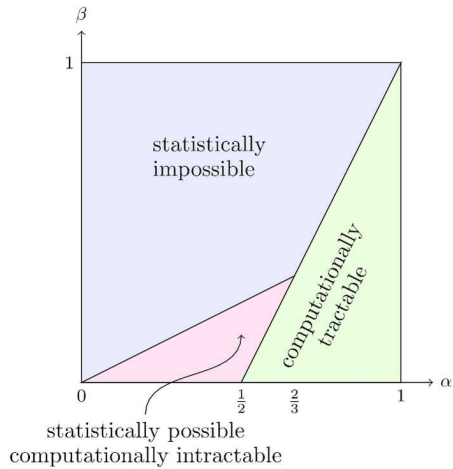}

\caption{Detection boundary $\beta^*$ versus efficiently computable
detection boundary $\beta^\sharp$.}
\label{figdiagram-test}
\end{figure}

Our main result can be illustrated by focusing on the following
asymptotic regime, where the submatrix size grows according to
$k=\Theta(p^{\alpha})$, and the signal magnitude decays as $\lambda=
\Theta(p^{-\beta})$ for fixed constants $\alpha\in(0,1)$ and $\beta
\in[0,1]$\footnote{The regime of $\beta>1$ is not interesting since
the hypotheses are indistinguishable even if the submatrix becomes the
whole matrix ($k=p$).} as $p\to\infty$.
For any two numbers $a$ and $b$, let $a\wedge b = \min(a,b)$ and
$a\vee b = \max(a,b)$.
Define
%
\[
\beta^* \triangleq\frac{\alpha}{2} \vee(2 \alpha- 1) \geq\beta^\sharp
\triangleq0 \vee(2 \alpha- 1).
\]
%
The statistical and computational feasibility of the submatrix
detection problem is demonstrated in Figure~\ref{figdiagram-test},
where the $(\alpha,\beta)$-plane is divided into three regions:
%
\begin{longlist}[(2)]
\item[(1)] $\beta> \beta^*$ (top region): reliable detection of the
submatrix is statistically impossible because the signal is too weak.
\item[(2)] $\beta< \beta^\sharp$ (right triangular region): reliable submatrix
detection is achievable by computationally efficient tests.
%
\item[(3)] $\beta^\sharp< \beta< \beta^*$ (lower left triangular region): reliable
detection is statistically possible but computationally intractable, in
the sense that it is at least as hard as solving the planted clique
problem of a particular configuration, which is intractable under the
planted clique hypothesis.
\end{longlist}
%
Therefore, the tractability of the submatrix detection problem
undergoes a sharp transition:
in the {moderately sparse} regime where $\alpha\in({2}/{3}, 1)$,
computational constraints incur no penalty on the statistical
performance. In contrast,
in the {highly sparse} regime where $\alpha\in(0, {2}/{3})$,
achieving the statistical optimal boundary requires computational
resources that are powerful enough to solve the planted clique problem,
and, consequently, computationally efficient procedures
require significantly higher signal-to-noise ratio to detect the submatrix.

The complexity-theoretic limits for the submatrix detection problem
also lead to interesting findings for the related support recovery
problem when the signal submatrix is present \cite{Shabalin09,BKRSW11,KBRS11}.
Moreover, it also sheds light on the statistical and computational
tradeoff in the problem of estimating block sparse matrices
\cite{MW13}.
In particular, we show the surprising result that the hardness of
minimax estimation can crucially depend on the loss function, in the
sense that attaining the minimax estimation rate can be computationally
easy for one type of loss functions but hard for the other.

\subsection{Related works}
Despite the vast body of literature on developing computationally
efficient procedures with optimal statistical performance for problems
such as compressed sensing, rigorous results on inferential limits of
statistical problems under computational complexity constraints are
comparatively limited.
A representative work is
{the investigation of the complexity of detecting sparse principal
components by Berthet and Rigollet \cite{Berthet13}, which is one of
the motivations of the present paper.
Sparse principal component detection refers to the problem of testing
$N(0,I_p)$ against $N(0,I_p+ a vv')$ for $k$-sparse unit vector $v$
based on $n$ i.i.d. observations \cite{Berthet12}.
Since the model is not discrete, as previously mentioned, the
difficulty of ill-defined computational complexity is also present.
In \cite{Berthet13}, the authors relaxed the Gaussian detection
problem to a composite testing problem that includes discrete
distributions, where the empirical projection variances of the null and
alternative hypotheses satisfy respective uniform $\chi^2$-tail type
concentration inequalities.
In the regime of $p^\delta\leq k \leq p^\alpha$ for some absolute
constants $0<\delta\leq\alpha< \frac{1}{2}$ and $n \leq p$, they
showed that the computable detection rate for the deviation of the
largest principal component from the rest of the spectrum is no smaller
than $\sqrt{\frac{k^b}{n}}$ for any $b < 2$, which far exceeds the
minimax detection rate of $\sqrt{\frac{k}{n}\log p}$.}

Although both the authors of current paper and Berthet and Rigollet
\cite{Berthet13} use the planted clique hypothesis for studying
complexity theoretic lower bounds, there are a few important differences.
First, Berthet and Rigollet \cite{Berthet13} extend the original
``simple vs. composite'' Gaussian sparse principal component detection
problem into a ``composite versus composite'' testing problem, and the
data no longer needs to be Gaussian.
As a consequence, more distributions are included in both the null and
alternative hypotheses, and thus constructing the reduction scheme
becomes easier than for the original Gaussian hypotheses.
In contrast, the current paper considers an asymptotically equivalent
discretized model which is faithful to the original Gaussian submatrix
detection problem in \cite{BI12}.
Second, the computational lower bounds in \cite{Berthet13} are
established only when the sparsity level satisfies $p^\delta\leq k\leq
p^\alpha$ for $0 < \delta\leq\alpha< \frac{1}{2}$.
In comparison, due to a new reduction scheme, the current paper
provides a more complete characterization of the computational limits
for all $k\geq p^\delta$ and any $\delta> 0$.
Last but not least, we propose an asymptotic equivalence framework in
the sense of Le Cam, which preserves the statistical nature of the
problem and, at the same time, allows rigorous statements of
computational complexity of testing procedures.
The approach via asymptotically equivalent discretized experiments is
potentially useful in future works dealing with nondiscrete distributions.

In addition, some researchers have studied the minimax sub-optimality
of certain computationally efficient methodologies, such as those based
on convex relaxations, in an array of problems including estimating
sparse eigenvectors \cite{KNV13}, support recovery for sparse
submatrices \cite{BKRSW11}, combinatorial testing \cite{Addario10},
community detection \cite{Arias13a}, etc.
In some of the papers, the authors also conjecture that the minimax
rate optimality cannot be achieved by any computationally efficient algorithms.
From a different viewpoint, Chandrasekaran and Jordan \cite{ChJo13}
consider the tradeoff between computation and statistical performance
within a specific family of
algorithms parameterized by the level of convex relaxations in the
classical normal mean estimation problem. In contrast, the goal of the
present paper is to investigate the impact of complexity constraint on
\emph{any} statistical procedure for the submatrix detection problem.

\subsection{Organization of the paper}

The rest of the paper is organized as follows.
In Section~\ref{secgaussian}, we study test statistics for submatrix
detection under Gaussian models.
To incorporate computational complexity into the decision theoretic
problem, we introduce in Section~\ref{secdiscrete} a sequence of
asymptotically equivalent discretized models and show that the minimax
detection results (\ref{eqminimax-possible})--(\ref
{eqminimax-impossible}) remain unchanged under these models.
In Section~\ref{seccomplexity}, we state our main result in
Theorem~\ref{thmad} under the planted clique hypothesis and present
its proof
with a concrete randomized polynomial-time procedure that reduces the
planted clique problem to a Bayesian version of the submatrix
detection problem.
We discuss some related problems in Section~\ref{secdiscuss}.
Section~\ref{secproof} presents additional proofs for results in
earlier sections.

\subsection{Notation}
For any positive integer $n$, let $[n]$ denote the set $ \{
1,\dots,n \}$.
For any $a\in\reals$, let $a_+ = a\vee0$.
For any square matrix $A$, $\tr(A) = \sum_{i} A_{ii}$ stands for its trace.
For any two matrices $A$ and $B$ of the same size, $A\circ B$ denotes
their component-wise product; that is, $(A\circ B)_{ij} = A_{ij} B_{ij}$,
and $\langle A, B\rangle= \tr(A'B)$.
Let $\calL(Y)$ denote the law, that is, the probability distribution,
of a random variable $Y$.
Let $\calL(Y\mid E)$ denote the distribution of $Y$ conditioned on the
event $E$.
The total variation distance between distributions ${P}$ and ${Q}$ is
$\TV(P,Q) \triangleq1-\int(\diff{P}\wedge\diff{Q})$.
For ease of notation, we also write $\TV(X,Y)$ in place of $\TV(\calL
(X), \calL(Y))$ for random variables $X$ and $Y$.
We write $X \stackrel{\mathrm{(d)}}{=} Y$ if $\calL(X) = \calL(Y)$.
Let $\Phi$, $\overline\Phi=1-\Phi$ and $\varphi$ denote the
distribution, survival and the probability density functions of the
standard Gaussian distribution.
For any set $I$, $\llvert I\rrvert$ denotes its cardinality.
For any sequences $\{a_p\}$ and $\{b_p\}$, we write $a_p\asymp b_p$ or
$a_p = \Theta(b_p)$ if there is an absolute constant $C>0$ such that
$1/C\leq a_p/b_p \leq C$.
We also write $a_p\ll b_p$ and $b_p\gg a_p$ if $a_p=o(b_p)$, and $a_p =
\Omega(b_p)$ if $b_p = O(a_p)$.

\section{Test statistics for submatrix detection}
\label{secgaussian}

To prepare for later investigation, we first study three test
statistics for the submatrix detection problem (\ref
{eqHgauss})--(\ref{eqset-H1+}).
The first two are
the linear and the scan test statistics proposed in \cite{BI12},
%
\begin{eqnarray}
\label{eqTlin-Tscan}
\Tlin &=& \Tlin(X) \triangleq\frac{1}{p} \sum
_{i,j=1}^p X_{ij},
\nonumber\\[-8pt]\\[-8pt]\nonumber
\Tscan &=& \Tscan(X)
\triangleq\frac{1}{k} \max_{\llvert S\rrvert = \llvert T\rrvert =
k} \sum
_{i \in S,j\in T} X_{ij}.
\end{eqnarray}
In addition, we also consider the maximum test statistic
%
\begin{equation}
\label{eqTmax} \Tmax= \Tmax(X) \triangleq\max_{i,j \in[p]}
X_{ij}.
\end{equation}

The following lemma gives nonasymptotic bounds on the Type-\textup{I${}+{}$II} error
probabilities on tests based on these statistics.
Recall the definition of $\calM(p,k,\lambda)$ in~(\ref{eqset-H1+}).

\begin{lemma}\label{lmmub}
Let $\calM= \calM(p,k,\lambda)$ and $c>0$ be any absolute constant.
For $\Tlin$ in (\ref{eqTlin-Tscan}),
set $\tau= \frac{\lambda k^2}{2p}$. Then
%
\begin{equation}
\Prob_0 \{ \Tlin> \tau\} + \sup_{\theta\in\calM}
\mathbb{P}_{\theta} \{ \Tlin\leq\tau\} \leq\eexp^{-\sfrac{\lambda^2
k^4}{8p^2}}.
\label{eqtlin}
\end{equation}
For $\Tscan$ in (\ref{eqTlin-Tscan}),
set $\tau'= \sqrt{ (4+c) \log{p\choose k}}$. Then
%
\begin{equation}
\Prob_0 \bigl\{ \Tscan> \tau' \bigr\} + \sup
_{\theta\in\calM} \mathbb{P}_{\theta} \bigl\{ \Tscan\leq
\tau' \bigr\} \leq\pmatrix{p
\cr
k}^{-c/2} +
\eexp^{- \sklfrac{1}{2} ( \lambda k -
\tau' )_+^2 }. 
\label{eqtscan}
\end{equation}
%
For $\Tmax$ in (\ref{eqTmax}),
set $\tau'' = \sqrt{(4+c)\log p}$. Then
%
\begin{equation}
\Prob_0 \bigl\{ \Tmax> \tau'' \bigr\} +
\sup_{\theta\in\calM} \mathbb{P}_{\theta} \bigl\{ \Tmax\leq
\tau'' \bigr\} \leq p^{-c/2} +
\eexp^{-\sklfrac{1}{2} ( \lambda- \tau'' )_+^2 }. \label{eqtmax}
\end{equation}
%
\end{lemma}

For the proof of Lemma~\ref{lmmub}, see Section~\ref
{secproof-test-statistics}.
By Lemma~\ref{lmmub}, we have\break $\calE({\mathbf{1}_{ \{{\Tlin>
\tau} \}}}) \to0$ when the first condition\vspace*{1pt} in (\ref
{eqminimax-possible}) holds,
while $\calE({\mathbf{1}_{ \{{\Tscan> \tau'} \}}})\to0$
when the second condition in (\ref{eqminimax-possible}) holds if we
pick the constant $c$ in $\tau'$ to be sufficiently small such that
$\liminf_{p\to\infty} \lambda k / \tau' > 1$.
The error bounds on $\Tmax$ will be used later to establish the
achievability part of the main result.

\section{Asymptotically equivalent discretized model}
\label{secdiscrete}

Gaussian distributions serve as good statistical models for many
real-world datasets.
However, as an idealized approximation, Gaussian experiment does not
capture the finite-precision nature of statistical computing systems in reality.
As mentioned in Section~\ref{secintro}, it is an ill-defined problem
to investigate the computational complexity of testing the Gaussian
hypothesis (\ref{eqHgauss}) since the data do not admit any
representation using finite bits.
Therefore a new paradigm is needed in order to
make sense of hypothesis testing with complexity constraints in general.
There are two goals of the paradigm:
%
\begin{longlist}[(a)]
\item[(a)] to provide a rigorous framework for quantifying the
complexity of statistical inference involving continuous, for example,
Gaussian, distributions and

\item[(b)] to preserve the statistical difficulty of the original
problem in the sense of Le Cam's asymptotic equivalence.
\end{longlist}
In this section, we propose such a paradigm based on discretizing the
original Gaussian experiment, which achieves both of the above goals.

\subsection*{Discretized models}

For any integer $t \in\naturals$, define the function $\quana
{t}{\cdot}\dvtx \reals\to2^{-t} \integers$ by
%
\begin{equation}
\quana{t} {x} = 2^{-t}{ \bigl\lfloor{2^t x} \bigr\rfloor}.
\label{eqquant}
\end{equation}
The function $\quana{t}{\cdot}$ naturally extends to matrices
componentwise: for $A = (A_{ij})$, $\quana{t}{A} = (\quana{t}{A_{ij}})$.

Recall the submatrix detection problem (\ref{eqHgauss}).
To model statistical inference with finite precision and complexity
constraints, let us consider the same testing problem based on the
discretized data $[X]_t$.
In other words, the hypotheses are
%
\begin{equation}
H_0^t\dvtx [X]_t \sim\Prob_0^{t}
\quad\mbox{versus}\quad H_1^t\dvtx [X]_t \sim
\Prob_{\theta}^{t},\qquad \theta\in\calM(p,k,\lambda), \label{eqHdis}
\end{equation}
where for $X\sim\Prob_\theta$,
\[
\Prob_\theta^{t} \triangleq\calL\bigl(\quana{t} {X}\bigr)
\]
is the discrete distribution induced by the quantization operation
(\ref{eqquant}), which is supported on $(2^{-t} \integers)^{p\times p}$.

Now on the discretized data, any test for (\ref{eqHdis}) is a
(possibly randomized) function from the countable set $(2^{-t}
\integers)^{p \times p}$ to $\{0,1\}$.
Since there exists a one-to-one mapping between the set $(2^{-t}
\integers)^{p \times p}$ and the set of all finite length binary
sequences $\bigcup_{n\in\naturals}\{0,1\}^n$, each observed $\tquana{t}{X}$
can be represented by a finite number of bits, and hence the
computational complexity of any test of interest is well defined; see,
for example, \cite{Arora09}, Chapter~7.
Thus the first goal of the paradigm is achieved.
As an aside, we note that though each coordinate of the discretized
data matrix $[X]_t$ has countably infinite support, its Shannon entropy
is finite and behaves according to $H([X_{ij}]_t) = t + O(1)$ as $t
\to\infty$ \cite{renyi}. Therefore, if we choose $t = \Theta(\log
p)$, then $[X]_t$ can be represented \emph{on average} using
variable-length lossless codes with $O(p^2 \log p)$ number of bits
\cite{cover}.

Given any test $\phi= \phi(\quana{t}{X})$ for (\ref{eqHdis}), we
can analogously define the worst-case Type-\textup{I${}+{}$II} error probability
$\calE(\phi)$ as in (\ref{eqperr-phi}) but with $\Prob_0$ and
$\Prob_\theta$ replaced by $\Prob_0^t$ and~$\Prob_\theta^t$.
Consequently, $\calE^*$ can also be defined as in (\ref{eqperr}).

\subsection*{Asymptotic equivalence}

Now we show that as long as we quantize each coordinate with accuracy
$p^{-c}$ for some constant $c>0$,
that is, $t = \Theta(\log p)$,
the resulting family of discretized distributions $\{\Prob_\theta^t\dvtx
\theta\in\reals^{p\times p}\}$ is asymptotically equivalent to the
original Gaussian experiment in the sense of Le Cam. Therefore any
inference problem, in particular, submatrix detection, performed on the
discretized data is asymptotically equally difficult as the original
problem as $p\to\infty$, and hence we also achieve the second goal of
the paradigm.

To state the equivalence result, recall the definition of Le Cam
distance between statistical experiments. Let $P$ be a probability
measure on a standard Borel space $(\sfX,\calF)$, and let $K$ denote
a probability transition kernel (Markov kernel) from $(\sfX,\calF)$
to a standard Borel space $(\sfY,\calG)$. Denote by $KP$ the
pushforward of $P$ under $K$, that is, $KP(\diff y) = \int_{\sfX}
K(\diff y\mid x) P(\diff x)$.
Given two experiments $\calP= \{P_{\theta}\dvtx \theta\in\Theta\}$ on
$(\sfX,\calF)$ and $\calQ= \{Q_{\theta}\dvtx \theta\in\Theta\}$ on
$(\sfY,\calG)$ with common parameter space $\Theta$, the \emph{Le
Cam deficiency} of $\calP$ with respect to $\calQ$ is defined by
\[
\delta(\calP,\calQ) \triangleq\inf_T \sup
_{\theta\in\Theta} \TV(T P_{\theta}, Q_{\theta}),
\]
where the infimum is over all probability transition kernels from
$(\sfX,\calF)$ to $(\sfY,\calG)$ \cite{SS00}, Theorem 1.7, page~29. The
\emph{Le Cam distance} between $\calP$ and $\calQ$ is
\[
\Delta(\calP,\calQ) \triangleq\delta(\calP,\calQ) \vee\delta(\calQ
,\calP).
\]
Two sequences of experiments $\{\calP^{(p)}\}_{p \in\naturals}$ and
$\{\calQ^{(p)}\}_{p \in\naturals}$ are \emph{asymptotically
equivalent} if their Le Cam distance vanishes \cite{LeCam86},
Section~2.3, that is, if
$\Delta(\calP^{(p)}, \calQ^{(p)}) \to0$ as $p \to\infty$.

The following theorem, proved in the supplement \cite{MW13-supp},
gives a nonasymptotic upper bound on the Le Cam distance between the
Gaussian experiments
$\calP^{(p)} = \{\Prob_\theta\dvtx \theta\in\reals^{p\times p}\}$
and its discretized version $\calP^{(p,t)} = \{\Prob_\theta^{t}\dvtx
\theta\in\reals^{p\times p}\}$. Therefore, as long as $t$ grows at a
logarithmical rate with $p$, the discretized model is asymptotically
equivalent to the original Gaussian model.

\begin{theorem}\label{thmequiv}
For any $t, p \in\naturals$,
$\Delta(\calP^{(p)}, \calP^{(p,t)}) \leq2 p^2 2^{-2 t/3}$.
Consequently, if $t=t(p) \geq(3+\varepsilon) \log_2 p$ for any
$\varepsilon> 0$, then $\{\calP^{(p)}\}_{p \in\naturals}$ and $\{\calP
^{(p,t(p))}\}_{p \in\naturals}$ are asymptotically equivalent as
$p\to\infty$.
\end{theorem}

For the proof of Theorem~\ref{thmequiv}, see the supplement \cite{MW13-supp}.
An immediate consequence of Theorem~\ref{thmequiv} on submatrix
detection is the following:
Since the difference between the optimal Type-\textup{I${}+{}$II} error probabilities
for the Gaussian hypotheses and the discretized hypotheses is upper
bounded by their Le Cam distance \cite{LeCam86}, Theorem 2.2, which
vanishes as $p \to\infty$, we conclude that testing on discretized data
has no impact on the statistical performance asymptotically in the
high-dimensional setting.
In particular, conclusion (\ref{eqminimax-possible})--(\ref
{eqminimax-impossible}) continues to hold for testing (\ref{eqHdis}).

%
\begin{remark}
\label{rmktest-discrete}
For the discretized model,
when either condition in (\ref{eqminimax-possible}) holds,
reliable detection can be attained by applying the linear or the scan
test to the quantized data.
To see this, note that the statistics $\Tlin, \Tscan$ and $\Tmax$
defined in (\ref{eqTlin-Tscan})--(\ref{eqTmax}) are all
$p$-Lipschitz with respect to the entrywise $\ell_\infty$-norm of $X$.
Using Lemma~\ref{lmmub},
it is straightforward to verify that if we compute $\Tlin$ and $\Tscan
$ based on the quantized data $\quana{t}{X}$ with $t \geq(3+\varepsilon
)\log p$,
then $\calE({\mathbf{1}_{ \{{\Tlin> \tau} \}}})$ [resp.,
$\calE({\mathbf{1}_{ \{{\Tscan> \tau'} \}}})$] vanishes
when the first (resp., second) condition in (\ref{eqminimax-possible}) holds.
Here, the thresholds $\tau$ and $\tau'$ are defined in Lemma~\ref{lmmub}.
\end{remark}

%
\begin{remark}
From an alternative viewpoint,
for appropriately chosen $t = t(p)\in\naturals$, one can restrict the
attention to all tests that are measurable with respect to the $\sigma
$-algebra on $\reals^{p\times p}$ generated\vspace*{1pt} by $\calF_t = \{\prod
_{i,j=1}^p [x_{ij}2^{-t}, (x_{ij}+1)2^{-t}), x_{ij}\in\integers\}$
rather than the usual Borel $\sigma$-algebra generated by all open sets.
Thus any such test $\psi$ remains constant on any set in $\calF_t$.
Moreover, $\psi(X) = \psi(\tquana{t}{X})$, and its computational
complexity is well defined.
Last but not least, the hypothesis testing problems (\ref{eqHgauss})
and (\ref{eqHdis}) become equivalent on this smaller $\sigma$-algebra.
\end{remark}

\section{Complexity theoretic limits}\label{seccomplexity}

In this section, we investigate complexity theoretic limits of the
submatrix detection problem by drawing its connection to the planted
clique problem. Let $N\in\naturals$ and $\kappa\in[N]$.
We denote by $\calG(N,1/2)$ the Erd\H{o}s--R\'enyi random graph on
$N$ vertices, where each edge is drawn independently at random with
probability $1/2$.
In addition, following \cite{Jerrum92,Alon98}, we use $\calG
(N,1/2,\kappa)$ to denote the random graph generated by first sampling
from $\calG(N,1/2)$, then picking $\kappa$ vertices
uniformly at random and connecting all edges in-between to form a
clique of size $\kappa$.
Distinguishing these two ensembles is known as the planted clique
problem, formally defined as follows:

\begin{definition}
\label{defpc}
Let $A \in\{0,1\}^{N\times N}$ be the adjacency matrix of a random
graph drawn from either $\calG(N,1/2)$ or $\calG(N,1/2,\kappa)$.
The \emph{planted clique problem of parameters $(N,\kappa)$},
denoted by $\PC(N,\kappa)$, refers to the hypothesis testing problem of
%
\begin{equation}
\label{eqpc-test} H_0^G\dvtx A\sim\calG(N,1/2)\quad\mbox{vs.}\quad H_1^G\dvtx A\sim\calG(N,1/2,\kappa).
\end{equation}
\end{definition}

The planted clique problem has a long history in the theoretical
computer science literature.
It is known that finding the clique is statistically impossible when
$\kappa= o(\log N)$. Moreover, a greedy algorithm succeeds if $\kappa
\geq c\sqrt{N \log N}$ for some constant $c>0$ \cite{Kuvcera95}.
Using spectral methods, Alon, Krivelevich and Sudakov \cite{Alon98}
provided the first polynomial time detection algorithm when $\kappa=
c\sqrt{N}$, with later improvements obtained in, for example, \cite
{Feige00,Feige10,Dekel10,Ames11,Deshpande13}.
However,
it is widely believed that the detection problem cannot be solved in
randomized polynomial time
when $\kappa= o(\sqrt{N})$,
which can be summarized as the following \emph{planted clique hypothesis}.
This version is similar to \cite{AAK07}, Conjecture 4.13, and \cite
{Berthet13}, Hypothesis $\mathsf B_{\PC}$.

\begin{hypo}
\label{hypopc}
For any sequence $\{\kappa_N\}$ such that $\limsup_{N \to\infty}
\frac{\log\kappa_N }{ \log N} < 1/2$
and any sequence of randomized polynomial-time tests\footnote
{Formally, randomized polynomial-time tests belong to the \textbf{BPP}
complexity class.
Interested readers are referred to standard textbooks on computational
complexity theory (e.g., \cite{Arora09}, Chapter~7) for the formal
definitions and discussions. Intuitively speaking, randomized
polynomial-time tests refer to algorithms with output space $ \{
0,1 \}$, which have access to external random numbers and whose running
time is bounded by a polynomial of the input length regardless of the
random numbers.}
$\{\psi_{N}\}$,
\begin{eqnarray*}
&&\liminf_{N\to\infty} \bigl( \Prob_{H_0^G} \bigl\{
\psi_{N}(A) = 1 \bigr\} + \Prob_{H_1^G} \bigl\{
\psi_{N}(A) = 0 \bigr\} \bigr) \geq\frac{2}{3}.
\end{eqnarray*}
\end{hypo}

Various hardness results in theoretical computer science have been
established based on the planted clique hypothesis, for example,
approximating the Nash equilibrium \cite{HK11}, independence testing
\cite{AAK07}, certifying the restricted
isometry property for compressed sensing measurement matrices \cite
{KZ12}, etc. Also, several cryptographic schemes have been proposed
assuming the intractability of finding planted cliques \cite
{Kuvcera92,JP00} or bicliques \cite{ABW10}.
Recently, the average-case hardness of planted clique has been established
under certain computation models; see, for example, \cite{Rossman10,Feldman13}.

The main result of the current paper is the following.

%
\begin{theorem}
\label{thmad}
Assume that Hypothesis~\ref{hypopc} holds.
Consider testing the discrete hypotheses
(\ref{eqHdis}) with $t = t(p) = 4{ \lceil{\log_2 p}
\rceil}$ in the asymptotic
regime (\ref{eqasymp}).
If, for some absolute constant $\delta> 0$,
%
\begin{equation}
\label{eqcomp-impossible} 
\frac{\lambda}{ p/k^{2+\delta}} \to0 \quad\mbox{and}\quad\limsup
_{p \to\infty} \lambda\sqrt{\log p} \leq\frac{1}{6},
\end{equation}
there exists \emph{no} sequence of randomized polynomial-time tests $\{
\phi_p\}$ such that $\liminf_{p \to\infty} \calE(\phi_p) < {2}/{3}$
for testing (\ref{eqHdis}).
Conversely, if
%
\begin{equation}
\label{eqcomp-possible} \frac{\lambda}{p/k^2 } \to\infty\quad\mbox
{or}\quad\liminf
_{p\to\infty} \frac{\lambda}{2\sqrt{\log p}} > 1,
\end{equation}
there is a sequence of linear-time tests $\{\phi_p\}$ such that $\calE
(\phi_p) \to0$.
\end{theorem}
%

As shown later in the proof of Theorem~\ref{thmad}, one can use
$\Tlin(\quana{t}{X})$ (resp., $\Tmax(\quana{t}{X})$) as the test
statistic when the first (resp., second) condition in (\ref
{eqcomp-possible}) holds.
It is straightforward to see that both $\Tlin$ and $\Tmax$ are of
linear complexity.

Contrasting the statistical limit (\ref{eqminimax-possible})--(\ref
{eqminimax-impossible}) with the computational limit (\ref
{eqcomp-impossible})--(\ref{eqcomp-possible}), we obtain the
following implication of Theorem~\ref{thmad} on the complexity of
submatrix detection:
suppose that $k \leq p^\alpha$ for some absolute constant $\alpha\in
(0,2/3)$. Then $\lambda\asymp\sqrt{\frac{1}{k}\log\frac{p}{k}}$
implies $\frac{\lambda}{p/k^2}\to0$. Consequently, conditions
(\ref{eqminimax-possible})--(\ref{eqminimax-impossible}) and
Theorem~\ref{thmequiv} imply that reliable detection is statistically
possible if and only if $\lambda= \Omega(\sqrt{\frac{1}{k}\log
\frac{p}{k}})$.
In contrast, condition (\ref{eqcomp-impossible}) in Theorem~\ref
{thmad} asserts that, to accomplish the same task using randomized
polynomial-time algorithms, it is necessary to have $\lambda= \Omega
(\frac{1}{\sqrt{\log p}} \wedge\frac{p}{ k^{2+\delta}})$ for all\vspace*{1pt}
$\delta> 0$, which far exceeds $\sqrt{\frac{1}{k}\log\frac{p}{k}}$
whenever $k \gg(\log p)^2$.
Therefore, computationally efficient procedures require significantly
larger signal level $\lambda$ to reliably detect the submatrix than
the statistical optimum.
More precisely, if $k = \Theta(p^\alpha)$ for some $\alpha\in
(0,2/3)$, then the minimal $\lambda$ for any randomized
polynomial-time test to succeed is at least $\lambda=\Omega(\frac
{1}{\sqrt{\log p}})$ when $\alpha\in(0,1/2)$ and $\Omega
(p^{1-2\alpha-\delta})$ for any $\delta> 0$ when $\alpha\in[1/2,
2/3)$, which
exceeds the statistical optimal level $\lambda= \Theta(p^{-\alpha
/2}\sqrt{\log p})$ by a polynomial factor in $p$.
Thus, in this regime, computational complexity constraints severely
limit the best possible statistical performance in the submatrix
detection problem.
On the other hand, when $k\geq p^\alpha$ for some $\alpha> {2}/{3}$,
$\frac{\lambda}{p/k^2}\to\infty$ is the dominating condition in
both (\ref{eqminimax-possible}) and (\ref{eqcomp-possible}), and a
computationally efficient test based on $\Tlin$ achieves statistically
optimal detection in this regime.
Figure~\ref{figdiagram-test} in Section~\ref{secintro} provides
a graphical illustration of the above discussion.

It should be noted that the sub-polynomial factor difference, that is,
${p}/{k^{2+\delta}}$ versus~${p}/{k^2}$, in the first part of (\ref
{eqcomp-impossible}) and (\ref{eqcomp-possible}) is a direct
consequence of Hypothesis~\ref{hypopc}.
In contrast, the logarithmic factor difference in the second part of
(\ref{eqcomp-impossible}) and (\ref{eqcomp-possible}) can
potentially be closed by employing better reduction argument and/or
more sophisticated testing procedures such as those based on spectral
methods, which we leave as a future direction.

The remainder of this section is devoted to proving Theorem~\ref
{thmad}, with auxiliary lemmas proved in Section~\ref{secproof}.
First, in Section~\ref{secreduce-one-sided} we provide some intuition
on how the planted clique problem is related to the submatrix
detection problem (\ref{eqHgauss}) under the Gaussian model.
Next, in Section~\ref{secdisc-reduction} we prove that under the
asymptotically equivalent discretized model,
every randomized polynomial time submatrix detector for (\ref
{eqHdis}) leads to a randomized polynomial time solver for the planted
clique problem of appropriate parameters with almost identical performance.
Finally, a proof of Theorem~\ref{thmad} is presented in Section~\ref
{seccomplexity-limit}.

\subsection{Planted clique and submatrix detection}\label{secreduce-one-sided}

We first explain how the submatrix detection problem can be reduced
from the planted clique problem under the original Gaussian model.
These results are presented as the precursor of the randomized
polynomial time reduction for the discretized model in Section~\ref
{secdisc-reduction}, as well as to provide insights into the hardness
of the submatrix detection problem.
A connection between the two problems has also been previously hinted
at in \cite{Addario10}.

Recall the Gaussian submatrix detection problem in (\ref{eqHgauss})
with parameter $(p,k,\lambda)$. For some $\ell\in\naturals$ to be
chosen depending on $p,k$ and $\lambda$, let
%
\begin{equation}
\label{eql} N = 2 p\ell.
\end{equation}
We construct a reduction scheme which maps any adjacency matrix $A\in\{
0,1\}^{N\times N}$ to a random matrix $X\in\reals^{p\times p}$ in
$O(N^2)$ number of flops, such that the following holds:
if $A$ is drawn from $\calG(N,1/2)$ under $H_0^G$, then the
distribution of $X$ is close in total variation distance to the null
distribution $\Prob_{0}$;
if $A$ is drawn from $\calG(N,1/2,\kappa)$ under $H_1^G$, then the
law of $X$ is close in total variation distance to a mixture of
distributions in the alternative $H_1$, where the clique size $\kappa$
is a constant multiple of $k$.


\subsubsection*{Randomized reduction}
An important step in
the following reduction scheme is to map any random edge to an $N(0,1)$
random variable and any edge in the clique to an $N(\mu,1)$ random
variable with some positive mean value $\mu$.
Although this goal might not be achievable exactly, we describe below a
strategy to achieve it approximately.

To this end,
for any $M \geq3$ and $0 < \mu\leq\frac{1}{2M}$, let $c_0 =
(1-2\overline\Phi(M))^{-1}$ and $c_1 = [1-\overline{\Phi}(M-\mu) -\overline
\Phi
(M+\mu) ]^{-1}$. We define two distributions $\calF_1$ and $\calF_0$
with the respective density functions
%
\begin{eqnarray}
\label{eqcond-gauss} f_1(x) &=& c_1 \varphi(x-\mu){
\mathbf{1}_{ \{{\llvert x\rrvert\leq M} \}
}},
\nonumber\\[-8pt]\\[-8pt]\nonumber
f_0(x) &=& \bigl[2c_0
\varphi(x) - c_1\varphi(x-\mu)\bigr]{\mathbf{1}_{ \{
{\llvert x\rrvert\leq M} \}}}.
\end{eqnarray}
Here both $f_0$ and $f_1$ are well-defined probability density functions.
In particular, the conditions $M\geq3$ and $0\leq\mu\leq\frac
{1}{2M}$ ensure that $f_0\geq0$.
In what follows, both $M$ and $\mu$, and thereby $\calF_0$ and $\calF
_1$, depend on $N$, though we suppress the dependence for notational
convenience.

The randomized mapping from $\bfA$ to $\bfX$ is as follows.
By (\ref{eql}), $N$ is even, and let $N_2 = N/2 = p\ell$ and
$[N]\setminus[N_2] = \{ N_2 + 1,\dots, N \}$:
%
\begin{longlist}[(2)]
\item[(1) (Gaussianization).]
Let $A_0 = A_{[N]\setminus[N_2], [N_2]}\in\reals^{N_2\times N_2}$
be the lower-left quarter of the matrix $A$.
Independent of $A$, we generate two $N_2\times N_2$ matrices $B_0$ and
$B_1$, whose entries are sampled independently from $\calF_0$ and
$\calF_1$ with density functions given in (\ref{eqcond-gauss}), respectively.
Define an $N_2\times N_2$ matrix $B$ by
%
\begin{equation}
B_{ij} = (B_0)_{ij} \bigl(1-(A_0)_{ij}
\bigr)+(B_1)_{ij} (A_0)_{ij}.
\label{eqAB01}
\end{equation}
In other words,
$\calL(B_{ij}\mid(A_0)_{ij} = 0) = \calF_0$ and $\calL(B_{ij}\mid
(A_0)_{ij} = 1) = \calF_1$.

\item[(2) (Partitioning).]
Partition $B$ into $\ell^2$ consecutive $p\times p$ blocks. In other
words, for $i,j\in[\ell]$, the $(i,j)$th block is $B^{(i,j)} =
(B^{(i,j)}_{a,b})\in\reals^{p\times p}$ where
%
\begin{equation}
\label{eqB-block} B^{(i,j)}_{a,b} = B_{(i-1)p+a, (j-1)p+b}\qquad\forall
a,b\in[p].
\end{equation}
%

\item[(3) (Averaging).]
Define $X\in\reals^{p\times p}$ by summing up all $\ell^2$ blocks
and scaling by~$\ell$:
%
\begin{equation}
\label{eqX-reduce} X = \frac{1}{\ell} \sum_{i= 1}^\ell
\sum_{j=1}^\ell B^{(i,j)}.
\end{equation}
\end{longlist}
Therefore, (\ref{eqAB01}), (\ref{eqB-block}) and (\ref
{eqX-reduce}) collectively define a deterministic function
%
\begin{eqnarray}\label{eqg}
g\dvtx \{0,1\}^{N \times N} \times\reals^{N_2 \times N_2} \times\reals
^{N_2 \times N_2} &\to& \reals^{p\times p},
\nonumber\\[-8pt]\\[-8pt]\nonumber
(A,B_0,B_1) &\mapsto& X
\end{eqnarray}
which can be computed in $O(N^2)$ number of flops.
The reason that we call the first step ``Gaussianization'' is due to
the following lemma, which ensures that for appropriately chosen $M$
and $\mu$, the marginal distribution of $B_{ij}$ is close to the
Gaussian distribution of unit variance and mean zero (resp., $\mu$) if
$(A_0)_{ij}$ corresponds to a random edge (resp., an edge in the clique).

\begin{lemma}
\label{lmmnormal-trunc}
Let $N\geq6$.
Let $\xi$ be a Bernoulli random variable. Let $W$ be a random variable
such that for $i\in\{0,1\}$, the conditional distribution of $W\mid
\xi=
i$ follows $f_i$ in (\ref{eqcond-gauss}) where $M\geq3$
and $\mu\leq\frac{1}{2M}$.
Then:
\begin{enumerate}
\item[(1)] if $\Prob\{ \xi= 1 \} = 1$, then $\TV(\calL
(W), N(\mu,1)) \leq\eexp^{(1-M^2)/2}$;
\item[(2)] if $\xi\sim\mathrm{Bernoulli}(1/2)$, then $\TV(\calL
(W), N(0,1)) \leq\eexp^{-M^2 /2}$.
\end{enumerate}
\end{lemma}
%

The following two lemmas characterize the law of $X=g(A,B_0, B_1)$ when
either $H_0^G$ or
$H_1^G$ in (\ref{eqpc-test}) holds.

%
\begin{lemma}
\label{lmmnull-reduce}
Suppose $H_0^G$ holds and $N\geq2p\geq6$.
Let $M\geq\sqrt{6\log N}$. Then
%
\begin{equation}
\TV\bigl(\calL(X), \Prob_0\bigr) \leq\frac{1}{p}.
\label{eqnull-reduce}
\end{equation}
\end{lemma}
%

%
\begin{lemma}
\label{lmmalter-reduce}
Suppose $H_1^G$ holds with $N \geq2p$, $p\geq2\kappa$ and $\kappa
\geq20$.
Let $k= { \lfloor{\kappa/20} \rfloor}$.
Let $M\geq\sqrt{6\log N}$ and $\mu\leq\frac{1}{2M}$ in (\ref
{eqcond-gauss}).
Then there exists a prior $\pi$ on $\calM= \calM(p,k, \frac{2\mu
p}{N})$ such that for $\bbP_\pi(\cdot) = \int_{\calM} \bbP_\theta
(\cdot)\pi(\diff\theta)$,
%
\begin{equation}
\TV\bigl(\calL(X), \bbP_\pi\bigr) \leq
\frac{1}{p}+
40k \biggl( \frac{\eexp}{4} \biggr)^{5k} + 2k \exp\biggl( -4k\log
\frac{p}{20k} \biggr). \label{eqalter-reduce}
\end{equation}
\end{lemma}
%

%
\begin{remark}
\label{rmkprior}
A careful examination of the proof of Lemma~\ref{lmmalter-reduce} in
Section~\ref{secproof-alter-reduce} reveals that the prior $\pi$ is
in fact supported on a subset $\widetilde\calM(p,k,\lambda) \subset
\calM
(p,k,\lambda)$ where
%
\begin{eqnarray}\label{eqset-prior}
\widetilde\calM(p,k,\lambda) &\triangleq& \bigl\{
\theta\in\reals^{p\times
p}\dvtx \exists U,V\subset[p],\mbox{ s.t. } k\leq
\llvert U\rrvert,\llvert V\rrvert\leq20k,
\nonumber\\[-8pt]\\[-8pt]\nonumber
&&\hspace*{5pt} \theta_{ij} \geq \lambda,\mbox{ if }(i,j)\in U\times V,
\theta_{ij} = 0\mbox{ if }(i,j)\notin U\times V \bigr\},
\end{eqnarray}
and $\lambda= \frac{2\mu p}{N}$.
In other words, any matrix in $\widetilde\calM(p,k,\lambda)$
contains a
nonzero rectangular submatrix whose row and column support sizes are
between $k$ and $20k$.
This observation will be useful for studying the hardness of estimation
in Section~\ref{secest-supp}.
\end{remark}

Combining Lemmas \ref{lmmnull-reduce} and \ref{lmmalter-reduce},
the following theorem shows that any submatrix detector leads to a test
with almost identical error probability for a planted clique problem,
whose parameters $(N,\kappa)$ depend on the parameters $(p,k,\lambda
)$ of the submatrix detection problem.

%
\begin{theorem}
Assume\vspace*{1pt} that $p\geq40k$ and $\lambda\leq\frac{1}{2\sqrt{6 \log(2p)}}$.
Suppose $\phi\dvtx \reals^{p\times p} \to\{0,1\}$ is a test for
distinguishing $H_0$ and $H_1$ in (\ref{eqHgauss})
with Type-\textup{I${}+{}$II} error probability upper bounded by $\varepsilon$, that is,
%
\begin{equation}
\mathbb{P}_{\theta_0} \bigl\{ \phi(X) = 1 \bigr\} + \sup
_{\theta
\in\calM(p,k,\lambda
)} \mathbb{P}_{\theta} \bigl\{ \phi(X) = 0 \bigr\}
\leq\varepsilon. \label{eqgood-test}
\end{equation}
%
Let $\kappa= 20 k$, $N = 2 p\ell$ and $N_2=N/2$, where $\ell$ is the
largest positive integer such that $N \sqrt{6 \log N} \leq
{p}/{\lambda}$.
Let $B_0,B_1\in\reals^{N_2\times N_2}$ have i.i.d. entries drawn from
$\calF_0$ and $\calF_1$, respectively, with $M = \sqrt{6\log N}$ and
$\mu= \frac{1}{2M}$.
Then $\psi(\cdot) = \phi(g(\cdot, B_0, B_1))$ is a test for the
planted clique detection problem (\ref{eqpc-test})
whose \mbox{Type-\textup{I${}+{}$II}} error probability is upper bounded by
%
\begin{equation}
\label{eqgbeta} \Prob_{H_0^G} \bigl\{ \psi(A) = 1 \bigr\} +
\Prob_{H_1^G} \bigl\{ \psi(A) = 0 \bigr\}\leq\varepsilon+ \beta,
\end{equation}
where $\beta= \frac{2}{p} + 40k ( \frac{\eexp}{4} )^{5k} +
2k\exp
( -4k\log\frac{p}{20k} )$.
\label{thmreduction}
\end{theorem}

\begin{pf}
Let $A$ denote the adjacency matrix of $\calG$.
By definition,
we have $M = \sqrt{6\log N}$, while the definition of $\ell$ ensures
that $\frac{2\mu p}{N}\geq\lambda$, and the constraint $\lambda\leq
\frac{1}{2\sqrt{6\log(2p)}}$ guarantees $\ell\geq1$.

By the definition of the total variation distance, Lemma~\ref
{lmmnull-reduce} implies that under~$H_0^G$,
%
\begin{eqnarray}
\label{eqbeta0} %
&& \bigl\llvert\Prob_{H_0^G}
\bigl\{ \phi\bigl(g(A,B_0,B_1)\bigr) = 1 \bigr\} - \Prob
_0 \bigl\{ \phi(X) = 1 \bigr\} \bigr\rrvert
\nonumber\\[-8pt]\\[-8pt]\nonumber
&&\qquad \leq\TV\bigl(\calL\bigl(g(A,B_0,B_1)\bigr),
\Prob_0\bigr) \leq\frac{1}{p} \triangleq\beta_0.
\end{eqnarray}
On the other hand, since $\frac{2\mu p}{N}\geq\lambda$, we have
$\calM(p,k,\frac{2\mu p}{N}) \subset\calM(p,k,\lambda)$.
So any location mixture $\Prob_\pi$ of the former can be viewed as a
mixture of the latter.
Hence, Lemma~\ref{lmmalter-reduce} implies that under $H_1^G$,
%
\begin{eqnarray}\label{eqbeta1}
&& \bigl\llvert\Prob_{H_1^G} \bigl\{ \phi\bigl(g(A,B_0,B_1)
\bigr) = 0 \bigr\} - \Prob_\pi\bigl\{ \phi(X) = 0 \bigr\} \bigr
\rrvert
\nonumber
\\
&&\qquad \leq \TV\bigl(\calL\bigl(g(A,B_0,B_1)\bigr),
\Prob_\pi\bigr)
\\
&&\qquad \leq
\frac{1}{p} + 40k \biggl( \frac{\eexp}{4} \biggr)^{5k} + 2k\exp
\biggl( -4k\log\frac{p}{20k} \biggr) \triangleq\beta_1.\nonumber
\end{eqnarray}
Since $\beta= \beta_0 + \beta_1$, the desired error bound (\ref
{eqgbeta}) follows from
\begin{eqnarray}
&& \mathbb{P}_{H_0^G} \bigl\{ \phi
\bigl(g(A,A_0,W)\bigr) = 1 \bigr\} + \mathbb{P}_{H_1^G} \bigl
\{ \phi\bigl(g(A,A_0,W)\bigr) = 0 \bigr\}
\nonumber
\\
&&\qquad \leq \mathbb{P}_{\theta_0} \bigl\{ \phi(X) = 1
\bigr\} + \mathbb{P}_{\pi} \bigl\{ \phi(X) = 1 \bigr\} +
\beta_0+\beta_1
\nonumber
\\
&&\qquad \leq\mathbb{P}_{\theta_0} \bigl\{ \phi(X) = 1 \bigr\} + \sup
_{\theta\in\calM
(p,k,\lambda)} \mathbb{P}_{\theta} \bigl\{ \phi(X) = 1 \bigr\} +
\beta
\nonumber
\\
&&\qquad \leq \varepsilon+ \beta,
\nonumber
\end{eqnarray}
where the last inequality is due to assumption (\ref{eqgood-test}) on
$\phi$.
\end{pf}

%
\begin{remark}
Although the reduction scheme $g$ can be implemented in $O(N^2)$ flops,
its computational complexity is ill defined
since it involves computing sums of continuous random variables and
processing infinitely many bits.
This issue will be addressed by a quantization argument in the next
subsection when we deal with the discretized models.
\end{remark}

\subsection{Randomized polynomial-time reduction for discretized models}\label{secdisc-reduction}

In this section, we show that with slight modifications, the scheme
introduced in Section~\ref{secreduce-one-sided} can be made into a
randomized polynomial-time reduction from the planted clique problem to
the submatrix detection problem under discretized models in rigorous
complexity-theoretic sense.

For the discretized model $\calP^{(p,t)}$ introduced in Section~\ref
{secdiscrete},
the reduction scheme from the planted clique model follows the same
steps in Section~\ref{secreduce-one-sided}, except that both the
input $(B_0,B_1)$ and the output $X$ are now discretized.

To this end, we first define discrete approximations, denoted by $Q_0$
and $Q_1$, to the densities $f_0$ and $f_1$ defined in (\ref{eqcond-gauss}).
Let $w,T$ be integers to be chosen based on $t,M$ and $N$. Recall the
quantization operator defined in (\ref{eqquant}) and that $B_0$ and
$B_1$ consist of i.i.d. entries drawn from densities $f_0$ and $f_1$,
respectively, which are supported on $[-M,M]$ by definition.
Note that each $\tquana{w}{(B_0)_{ij}}$ is drawn from a distribution
with atoms $x_i$ and probability mass function (p.m.f.) $p_i$ for $i\in[M
2^{w+1}]$.
To find a dyadic\vspace*{1pt} approximation for the p.m.f., let $q_i = \lfloor p_i 2^T
\rfloor2^{-T}$ for $i=2, \ldots, M2^w$ and $q_1 = \lfloor p_1
2^T\rfloor2^{-T} + 1 - \sum_{i \geq2} q_i$, where
%
\begin{equation}
T = { \lceil{\log_2 M} \rceil}+w + 3 \log_2 N.
\label{eqT}
\end{equation}
Denote by $Q_0$ the discrete distribution with atoms $x_i$ and
probability masses $q_i$. Similarly, define $Q_1$ as the dyadic
approximation for the distribution of $\tquana{w}{(B_1)_{ij}}$.

The reduction scheme operates as follows: first, generate $\breve
{B_i}$ consisting of i.i.d. entries drawn from $Q_i$ for $i=0,1$.
Next, replace the matrices $B_0$ and $B_1$ in (\ref{eqAB01}) by their
discretized version $\breve{B}_0$ and $\breve{B}_1$, and denote the
resulting matrix by $\breve B$.
Applying (\ref{eqB-block})--(\ref{eqX-reduce}) to $\breve B$, we
obtain $\breve{X}$ and output its quantized version $\tquana{t}{\breve{X}}$.
Implementing the above steps yields a deterministic function
%
\begin{eqnarray}\label{eqgbreve}
\qquad\breve{g}\dvtx \{0,1\}^{N \times N} \times\bigl([-M, M]_w
\bigr)^{N_2 \times N_2} \times\bigl([-M, M]_w\bigr)^{N_2 \times N_2} &\to&
\bigl(2^{-t}\integers\bigr)^{p\times p},
\nonumber\\[-8pt]\\[-8pt]\nonumber
(A,\breve B_0,\breve B_1) &\mapsto& \tquana{t} {
\breve{X}},
\end{eqnarray}
where $[-M, M]_w = [-M, M]\cap2^{-w}\integers$ is the quantized interval.

%
\begin{remark}[(Computational complexity of reduction)]\label
{rmkreduction-discrete}
First we discuss the complexity for generating the auxiliary random
variables used in the reduction scheme.
Note that each $(\breve B_0)_{ij}$ is drawn from $Q_0$ whose atoms
$x_i$ can be represented by ${ \lceil{\log_2 M} \rceil
}+w$ bits and the p.m.f.
$q_i$ is a dyadic rational with $T$ bits. Therefore sampling from the
distribution $Q_0$ can be done using the inverse CDF\footnote{More
sophisticated random number generators for discrete distributions (such
as Walker's alias method which requires linear time for preprocessing
and constant time per sample) can be found in~\cite{taocp2},
Section~3.4.1.} by outputting $x_J$, where $J = \min\{j\dvtx
\sum_{i=1}^j q_i \leq U 2^{-T} \}$ and $U$ is a random integer
uniformly distributed on $[2^T]$. Consequently, sampling from $Q_0$
requires $O(M 2^w T)$ preprocessing time to compute the CDF, and $T$
fair coin flips and $O(\log M + w)$ time per sample (via binary search).
Furthermore, discretizing each entry $\breve{X}_{ij}$ to $\tquana
{t}{\breve{X}_{ij}}$ involves keeping the first $t$ bits after the
binary point, which can be computed in $O(t)$ time.
Therefore we conclude that $\breve{g}$ can be computed using
$O(({ \lceil{\log_2 M} \rceil}+w+t)N^2)$ number of
binary operations.

To summarize, the randomized reduction scheme requires $O(N^2 T)$
random bits and $O(M 2^w T + N^2 (\log M + w+t))$ computation, where
$T$ is defined in (\ref{eqT}).
As we will show in Section~\ref{seccomplexity-limit}, for all cases
of interest in this paper,
we can set
$N = O(p^{2})$
and $M,w,t = O(\log_2{p}) = O(\log_2 {N})$.
Therefore,
our reduction for discretized models $A \mapsto\breve{g}(A, \breve
B_0, \breve B_1)$ is a \emph{randomized polynomial-time reduction}.
\end{remark}

We now investigate the distributions of $\tquana{t}{\breve{X}}$ under
$H_0^G$ and $H_1^G$, respectively.
The following lemmas are counterparts of Lemmas \ref{lmmnull-reduce}
and \ref{lmmalter-reduce} for discretized models.
Comparing with the total variation bound (\ref{eqnull-reduce}) and
(\ref{eqalter-reduce}), we show that, upon suitable choices of $w$
depending on $(t,N)$, replacing $B_0$ and $B_1$ with the discrete
versions $\breve{B}_0$ and $\breve{B}_1$ only introduces an extra
term of $4/p$ in the total variation
of $\calL(\tquana{t}{\breve{X}})$ to $\bbP_0^t$ under $H_0^G$, and
to a mixture of the alternative distributions $\bbP_\theta^t$ under
$H_1^G$, respectively.
This objective is accomplished by putting the support of $Q_0$ and
$Q_1$ on a finer grid than that of the output $\bX$, that is, choosing
$w > t$, which is essential for controlling the approximation error in
the output distribution incurred by quantizing the input.

%
\begin{lemma}
\label{lmmnull-discrete}
Let $N\geq2p\geq6$.
Let $w\in\naturals$ satisfy
%
\begin{eqnarray}\label{eqmw-alter}
w & \geq& t + 6 \log_2 N.
\end{eqnarray}
Then under $H_0^G$,
%
\begin{equation}
\TV\bigl(\calL\bigl(\tquana{t} {\breve{X}}\bigr), \bbP_0^t\bigr)
\leq\frac{5}{p}. \label{eqnull-discrete}
\end{equation}
\end{lemma}

%
\begin{lemma}
\label{lmmalter-discrete}
Suppose $H_1^G$ holds with $N \geq2p$, $p\geq2\kappa$, $\kappa\geq
20$, $k= { \lfloor{\kappa/20} \rfloor}$.
Let $M\geq\sqrt{6\log N}$ and $\mu\leq\frac{1}{2M}$. Let $w$
satisfy (\ref{eqmw-alter}). Then there exists a prior $\pi$ on
$\calM(p,k,\frac{2 \mu p}{N})$, such that for $\bbP_\pi^t(\cdot) =
\int_{\calM}\bbP_{\theta}^t(\cdot) \pi(\diff\theta)$,
%
\begin{eqnarray}
\TV\bigl(\calL\bigl(\tquana{t} {\breve{X}}\bigr), \bbP_\pi^t\bigr)
&\leq& \frac{5}{p}+ 40k \biggl( \frac{\eexp}{4} \biggr)^{5k} +
2k \exp\biggl( -4k\log\frac
{p}{20k} \biggr). \label{eqalternative-discrete}
\end{eqnarray}
\end{lemma}

Combining the two lemmas, we obtain the following result analogously to
Theorem~\ref{thmreduction}.

%
\begin{theorem}
Assume that $p\geq40k$ and $\lambda\leq\frac{1}{2\sqrt{6 \log(2p)}}$.
Suppose
$\phi\dvtx\break  (2^{-t}\integers)^{p\times p} \to\{0,1\}$ is a test for
distinguishing $H_0^t$ and $H_1^t$ in (\ref{eqHdis})
with Type-\textup{I${}+{}$II} error probability upper bounded by $\varepsilon$, that is,
%
\begin{equation}
\Prob_{0}^t \bigl\{ \phi\bigl(\tquana{t} {X}\bigr) = 1 \bigr\} +
\sup_{\theta
\in\calM
(p,k,\lambda)} \Prob_\theta^t \bigl\{ \phi\bigl(
\tquana{t} {X}\bigr) = 0 \bigr\} \leq\varepsilon. \label{eqgood-test-discrete}
\end{equation}
Let $\kappa,N,N_2\in\naturals$
be chosen as in Theorem~\ref{thmreduction}.
Let $w\in\naturals$ satisfy (\ref{eqmw-alter}) and $\breve{g}\dvtx\{
0,1\}^{N \times N} \times([-M, M]_w)^{N_2 \times N_2} \times([-M,
M]_w)^{N_2 \times N_2} \to(2^{-t}\integers)^{p\times p}$ be
defined in (\ref{eqgbreve}).
Then $\psi(\cdot) = \phi(\breve{g}(\cdot, \breve{B}_0, \breve
{B}_1))$ is a test for the planted clique detection problem (\ref{eqpc-test})
whose Type-\textup{I${}+{}$II} error probability is upper bounded by
%
\begin{equation}
\label{eqgbeta-discrete} \Prob_{H_0^G} \bigl\{ \psi(A) = 1 \bigr\} +
\Prob_{H_1^G} \bigl\{ \psi(A) = 0 \bigr\}\leq\varepsilon+ \beta,
\end{equation}
where $\beta= \frac{10}{p} + 40k ( \frac{\eexp}{4} )^{5k} +
2k\exp
( -4k\log\frac{p}{20k} )$.
\label{thmreduction-discrete}
\end{theorem}

\begin{pf}
In view of the analogy between Lemmas \ref{lmmnull-reduce}--\ref{lmmalter-reduce} and \ref{lmmnull-discrete}--\ref{lmmalter-discrete},
the proof follows the same argument as that in the proof of
Theorem~\ref{thmreduction}, except that $\Prob_0$ and $\Prob_\pi$
are replaced by $\Prob_0^t$ and $\Prob_\pi^t$, respectively, $g(A,
B_0, B_1)$ is replaced by $\breve{g}(A, \breve{B}_0, \breve{B}_1)$,
$\beta_0$ and $\beta_1$ in (\ref{eqbeta0}) and (\ref{eqbeta1})
are both increased by $4/p$.
\end{pf}

\subsection{Proof of Theorem~\texorpdfstring{\protect\ref{thmad}}{2}}
\label{seccomplexity-limit}

We start with the lower bound.
Without loss of generality, we can assume that $\lambda\geq{1}/{p}$
since when $\lambda< {1}/{p}$, both conditions in (\ref
{eqminimax-impossible}) hold in the asymptotic regime (\ref
{eqasymp}), and the problem is statistically impossible.
Let the sequence $\{(k(p),\lambda(p))\}$ satisfy (\ref{eqasymp}) and
(\ref{eqcomp-impossible}).
Let $\{\phi_{p}\}$ be a sequence of randomized polynomial time tests.
For conciseness we drop the indices in $k(p),\lambda(p)$ and $\phi_p$.
Suppose for the sake of contradiction that
%
\begin{equation}
\label{eqcor-compleixty-conclusion} \liminf_{p\to\infty} \Bigl(\Prob_{0}^t
\bigl\{ \phi\bigl(\tquana{t} {X}\bigr) = 1 \bigr\} + \sup_{\theta\in\calM
(p,k,\lambda)}\Prob
_\theta^t \bigl\{ \phi\bigl(\tquana{t} {X}\bigr) = 0 \bigr\} \Bigr) <
\frac{2}{3}. 
\end{equation}
Choose $\kappa$ and $N$ as in Theorems \ref{thmreduction}--\ref
{thmreduction-discrete}, that is, $\kappa= 20k$ and $N=2p\ell$ where
$\ell$ is the largest integer such that $ N \sqrt{6 \log N} \leq
{p}/{\lambda}$.
Since the first condition in~(\ref{eqcomp-impossible}) implies that
$\lambda\leq{Cp}/{k^{2+\delta}}$ for some constant $C$ and all
sufficiently large $p$, we have
$2 \kappa^{2+\delta/2} \sqrt{6 \log(2 \kappa^{2+\delta/2})} \leq
{p}/{\lambda}$ and hence $\ell\geq{ \lfloor{{\kappa
^{2+\delta/2}}/{p}} \rfloor}$ for all sufficiently large $p$.
Similarly, the second
condition in (\ref{eqcomp-impossible}) implies that $\lambda\leq
\frac{1}{2 \sqrt{6\log(2p)}}$ for all sufficiently large $p$ and
consequently, $\ell\geq1$.
Using the simple fact that $x { \lfloor{{y}/{x}} \rfloor}
\vee x \geq{y}/{2}$ for
all $x,y>0$, we conclude that $N = 2 p \ell\geq\kappa^{2+\delta/2}
\vee(2p)$
for all sufficiently large $p$, hence
%
\begin{equation}
\liminf_{p \to\infty} \frac{\log\kappa}{\log N} \leq\frac
{1}{2+\delta/2} <
\frac{1}{2}. \label{eqPCad}
\end{equation}
On the other hand, we have $N\leq{p}/{\lambda} \leq p^2$, where the
last inequality holds since we have assumed $\lambda\geq{1}/{p}$.
Applying Theorem~\ref{thmreduction-discrete} with $w= 16{ \lceil
{\log_2 p} \rceil} \geq t + 12\log_2 p \geq t + 6\log_2 N$,
we conclude from
(\ref{eqcor-compleixty-conclusion}) that the randomized test $\psi
(\cdot) = \phi(\breve{g}(\cdot, \breve{B}_0, \breve{B}_1))$ satisfies
%
\begin{equation}
\liminf_{p\to\infty} \bigl( \Prob_{H_0^G} \bigl\{ \psi(A) = 1
\bigr\} + \Prob_{H_1^G} \bigl\{ \psi(A) = 0 \bigr\} \bigr) <
\frac{2}{3}. \label{eqpsi}
\end{equation}
In view of Remark~\ref{rmkreduction-discrete}, $A\mapsto\breve
{g}(A, \breve{B}_0, \breve{B}_1)$ is a randomized polynomial-time reduction.
By the assumption on $\phi$, $\psi$ as a composition of $\phi$ and
$\breve{g}$ is a randomized polynomial-time test for $\PC(N,\kappa)$.
Therefore, (\ref{eqpsi}) contradicts Hypothesis~\ref{hypopc} in
view of (\ref{eqPCad}).

It remains to show the upper bound.
Denote the linear and maximum test statistics computed on the
discretized matrix $[X]_t$ by $\Tlin$ and $\Tmax$, respectively.
If the first condition in (\ref{eqcomp-possible}) holds, that is,
$\lambda k^2/p\to\infty$, in view of Lemma~\ref{lmmub} and
Remark~\ref{rmktest-discrete}, we have $\calE({\mathbf{1}_{ \{
{\Tlin> \tau} \}}}) \to0$ where $\tau$ is defined in
Lemma~\ref{lmmub}.
If the second condition in (\ref{eqcomp-possible}) holds,
recall $\tau'' = \sqrt{(4+c)\log p}$ defined in Lemma~\ref{lmmub}.
If the constant $c$ is sufficiently small such that $\liminf_{p\to
\infty} \lambda/ \tau'' > 1$, then following the reasoning in
Remark~\ref{rmktest-discrete},
it is straightforward to verify that $\calE({\mathbf{1}_{ \{
{\Tmax> \tau''} \}}}) \to0$.
This completes the proof.

\section{Discussion}
\label{secdiscuss}

In this paper, assuming the planted clique hypothesis, we have
demonstrated a phase transition phenomenon on gaps between
the optimal statistical performance with and without computational
complexity constraints
for the submatrix detection problem.
The hardness result in Theorem~\ref{thmad} has important consequences on
the hardness of two related problems, namely, \emph{support recovery}
and \emph{matrix estimation} under submatrix sparsity, both of which
are more difficult than detection and require stronger signal level.
To discuss computational complexity of statistical procedures, we focus
on the discretized models introduced in Section~\ref{secdiscrete}
throughout the current section.

\subsection{Support recovery}
\label{secsupp}

As previously studied in \cite{KBRS11}, the goal of support recovery is
to identify the minimum $\lambda$ such that, under the alternative
hypothesis $H_1$ in (\ref{eqHgauss}), the submatrix can be
consistently located.
According to Theorems 1 and 2 of \cite{KBRS11}, for all $k\leq p/2$,
one needs $\lambda= \Omega(\sqrt{ \log(p)/k})$
to recover the support consistently under the parameter space (\ref
{eqset-H1+}).
Compared with (\ref{eqminimax-possible})--(\ref
{eqminimax-impossible}), this coincides with the minimum signal
strength required for detecting the submatrix when $k = O(p^\alpha)$
for $\alpha< {2/3}$, but is much larger when $k = \Omega(p^\alpha)$
for $\alpha> {2/3}$.

Intuitively, locating the submatrix is more difficult than detecting
the mere existence thereof.
Therefore, the complexity theoretic limit for support recovery should
also exceed that of detection.
This claim, however, does not follow immediately since support recovery
only deals with the alternative hypothesis (\ref{eqset-H1+}), and the
null hypothesis is excluded from the parameter set.
To provide a rigorous argument, for $\lambda= \Omega(\sqrt{\log
(p)/k})$, given a support estimator $(\widehat{U}, \widehat{V})$ such that
$\sup_{\theta\in\calM} \Prob_\theta\{(\widehat{U}, \widehat{V}) \neq
(U,V)\} \leq\varepsilon$, we can construct a test for (\ref{eqHgauss})
which rejects if $T = \sum_{i\in\widehat{U}, j\in\widehat{V}} X_{ij}$
exceeds $\tau'$ defined in Lemma~\ref{lmmub}. Since $T$ is at most
$\Tscan$,
the same argument in the proof of Lemma~\ref{lmmub} shows that the
Type-\textup{I${}+{}$II} error probability for this test is upper bounded by
$\varepsilon
$ plus the right-hand side of (\ref{eqtscan}), which vanishes when
$p,k \to\infty$.
This implies that the minimal $\lambda$ achievable by computationally
efficient support estimators is at least a constant factor of that
required by computationally efficient submatrix detectors.
Therefore, Theorem~\ref{thmad} implies that no randomized
polynomial-time algorithm can achieve consistent support recovery when
condition (\ref{eqcomp-impossible}) holds.
This resolves in the negative the open question raised in \cite
{BKRSW11}, Section~5, on the existence of computationally efficient
minimax procedures in the regime of $k = O(p^\alpha)$ for any $\alpha
< {2}/{3}$.
It remains open to determine whether the statistically optimal support
recovery can be achieved computationally efficiently when $k = \Theta
(p^\alpha)$ for $\alpha> {2}/{3}$.

\subsection{Hardness of estimation depends on norm}
\label{secest-supp}

We now consider the computational aspect of the related problem of
estimating the mean matrix with submatrix sparsity under squared norm
losses. Denote the set of $k\times k$-sparse matrices by
%
\begin{eqnarray}\label{eqpara-space}
\calF(p,k) &=& \bigl\{ \theta\in
\RR^{p\times p}\dvtx  \exists U,V\subset[p],\mbox{ s.t.
}\llvert U\rrvert,
\llvert V\rrvert\leq k,
\nonumber\\[-8pt]\\[-8pt]\nonumber
&&\hspace*{83pt} \theta_{ij} = 0\mbox{ if } (i,j)\notin U\times V\bigr\},
\end{eqnarray}
which includes both the zero matrix and the set $\widetilde\calM(p,
{ \lfloor{k/20} \rfloor}, \lambda)$ [defined in (\ref
{eqset-prior})] for any $\lambda
> 0$.

Given the noisy observation $X = \theta+ Z$, where $Z$ consists of
standard Gaussian entries, the minimax risk
%
\begin{equation}
\label{eqrate} \Psi_{\llVert\cdot\rrVert}(p,k) \triangleq\inf_{\tilde
{\theta}}
\sup_{\theta\in\calF(p,k)} \Expect\llVert\tilde{\theta}-\theta\rrVert
^2
\end{equation}
has been obtained in \cite{MW13}, Section~4, within universal constant
factors using convex geometry and information-theoretic arguments for
all unitarily invariant norms,\footnote{To be precise, note that the
minimax rates in \cite{MW13} are obtained for the Gaussian model.
Since the loss function is unbounded, one cannot directly conclude from
asymptotic equivalence that the same rate applies to the discretized
model. Nevertheless, it is straightforward to extend the arguments in~\cite{MW13}, Section~4.1, to show that the rate of $\Psi_{\llVert
\cdot\rrVert
}(p,k)$ applies to the discretized model in Section~\ref{secdiscrete}
as long as $t = \Omega(\sqrt{\log p})$, independent of the unitarily
invariant norm. In particular, the lower bound in \cite{MW13},
Section~4.1.1, applies verbatim due to the data processing
inequality of the KL divergence, which is attained by the same
estimator defined in \cite{MW13}, Section~4.1.2, if $2^{-t} \leq
1/\sqrt{p}$.}
in particular, satisfies
%
\[
k \llVert I_k\rrVert^2 \lesssim\Psi_{\llVert\cdot
\rrVert}(p,k)
\lesssim k \llVert I_k\rrVert^2 \log
\frac{\eexp p}{k}.
\]
%

Capitalizing on the hardness result of detecting submatrices in
Theorem~\ref{thmad}, we show that the minimax estimation rates
corresponding to certain norm losses cannot be attained by
computationally efficient methods.
For conciseness, let us focus on the class of Schatten-$q$ norms
$\llVert\cdot\rrVert_{\mathrm{S}_q}$, defined as the $\ell
_q$-norm of singular values
for $q\in[1,\infty]$. The minimax rate is given by (see \cite{MW13},
Example 1)
%
\begin{equation}
\label{eqreg-rate-Sch} \Psi_{\mathrm{S}_q}(p,k) \asymp k^{{2}/{q}+1} +
k^{(2/q)\vee1 } \log\frac{\eexp p}{k},
\end{equation}
which is within a logarithmic term of $k^{{2}/{q}+1}$.
Next we discuss the computational cost of estimation by focusing on the
asymptotic regime where $k = \Theta(p^\alpha)$ for some $\alpha\in(0,1)$.
In view of the relationship between testing and estimation, we can use
the construction in Lemmas \ref{lmmnull-reduce}--\ref
{lmmalter-reduce} for the detection problem
\[
H_0\dvtx \theta= 0 \quad\mbox{versus} \quad H_1\dvtx \theta\in
\widetilde\calM\bigl(p, { \lfloor{k/20} \rfloor}, \lambda\bigr)
\]
as a two-point lower bound.
Note that for any $\theta\in\widetilde\calM(p, { \lfloor
{k/20} \rfloor},
\lambda)$
and any $q\in[1,\infty]$, $\llVert\theta\rrVert_{\mathrm
{S}_q} \geq\llVert\theta\rrVert_{\mathrm{S}_\infty} = \Omega(k\lambda)$.
Assuming Hypothesis~\ref{hypopc}, we conclude that the squared
Schatten-$q$ norm risk achievable by any randomized polynomial-time
estimator is at least $\Omega(k^2 \lambda^2)$ for any $\lambda$
satisfying (\ref{eqcomp-impossible}).
Thus, for any constant $\delta> 0$, the worst-case risk is at least
$\Omega(k^{-\delta}(k^2 \wedge\frac{p^2}{k^2}))$.
Note that this lower bound is not monotonic in $k$ and can be easily
improved to
%
\begin{equation}
\Omega\bigl(k^{-\delta}\bigl(k^2 \wedge p\bigr) \bigr)
\label{eqestimation-comp1}
\end{equation}
since the risk is clearly nondecreasing in $k$.

On the constructive side, an estimation error of
%
\begin{equation}
O \biggl( k^{(2/q+1)\vee2} \log\frac{\eexp p}{k} \wedge p^{2/q+1}
\biggr) \label{eqestimation-comp2}
\end{equation}
in squared Schatten-$q$ norm can be achieved in polynomial time. To see
this, first note that treating a $k \times k$-sparse matrix as a
$k^2$-sparse vector in $p^2$-dimensional space and applying entrywise
hard thresholding yields an estimator $\hat{\theta}$ whose
mean-square error (i.e., squared Frobenius or Schatten-2 norm) is at
most $O(k^2 \log\frac{\eexp p}{k})$. Then we project $\hat{\theta}$
into the space of row-sparse matrices to obtain $\tilde{\theta}$ by
choosing the $k$ rows of $\hat{\theta}$ of the largest $\ell_2$-norm
and set the remaining rows to zero. Since the estimand $\theta$ also
has $k$ nonzero rows,\vspace*{1pt} applying the triangle inequality yields\vadjust{\goodbreak} $\Expect
\llVert\tilde{\theta} -\theta\rrVert_{\mathrm{S}_2}^2 =
O(k^2 \log\frac{\eexp
p}{k})$, which implies $\Expect\llVert\tilde{\theta} -\theta
\rrVert_{\mathrm{S}_q}^2 = O(k^{(2/q+1)\vee2} \log\frac{\eexp p}{k})$, since $\llVert
\cdot
\rrVert_{\mathrm{S}_q} \leq(1 \vee k^{1/q-1/2}) \llVert\cdot
\rrVert_{\mathrm{S}_2}$ for
all rank-$k$ matrices. Finally, simply estimating $\theta$ by the
observation $X$ achieves $O(p^{2/q+1})$.

Comparing the minimax rate (\ref{eqreg-rate-Sch}) with the
computationally lower and upper bounds (\ref
{eqestimation-comp1})--(\ref{eqestimation-comp2}), we obtain the
following result, assuming Hypothesis~\ref{hypopc}:
\begin{itemize}
\item For $q \in[1,2]$, using the entrywise thresholding estimator
defined above, the minimax rate is attained within a logarithmic factor
simultaneously for all $k$;
\item For $q \in(2,\infty]$, the minimax rate $\Psi_{\mathrm{S}_q}(p,k)$ cannot be attained by computationally efficient estimator
if $k = \Theta(p^\alpha)$ for all $\alpha\in(0, {\frac{q}{2+q}})$.
In this regime, entrywise thresholding is optimal within a
sub-polynomial factor among all randomized polynomial-time procedures.
\end{itemize}
%
More generally, one can show that for all \emph{quadratic norms} (see
\cite{Bhatia}, page~95), entrywise thresholding is near optimal
(within a sub-polynomial factor) among all computationally efficient
estimators. This extends the above result since Schatten\mbox{-$q$} norm is
quadratic if and only if $q \in[2,\infty]$.

\section{Proofs}
\label{secproof}
We present below the proofs of Lemmas \ref{lmmub}--\ref
{lmmalter-reduce}. The proofs of Theorem~\ref{thmequiv} and Lemmas
\ref{lmmnull-discrete} and \ref{lmmalter-discrete} are deferred to
the supplement \cite{MW13-supp}.
\subsection{Proof of Lemma~\texorpdfstring{\protect\ref{lmmub}}{1}}
\label{secproof-test-statistics}
Under $\Prob_0$, $\Tlin\sim N(0,1)$, hence $\Prob_0 \{ \Tlin>
\tau\} = \overline{\Phi}(\tau)$.
Under $\Prob_\theta$ for any $\theta\in\calM$, $\Tlin\sim N(\bar\theta
,1)$, where $\bar\theta\triangleq\frac{1}{p}\sum\theta
_{ij} \geq\frac{k^2 \lambda}{p}$
by the definition of $\calM$. Therefore $\Prob_{\theta} \{
\Tlin\leq\tau\} \leq\overline{\Phi}(\tau)$. Then (\ref
{eqtlin}) follows in
view of the Chernoff bound $\overline{\Phi}(\tau) \leq\frac{1}{2}\exp
(-\tau^2/2)$.

By the union bound,
$\Prob_0\{\Tscan> \tau'\}
\leq\binom{p}{k}^2 \Prob_0\{\sum_{i,j=1}^k{X_{ij}} > k\tau'\}
\leq \binom{p}{k}^2\times\break   \exp(-\tau'^2/2)
\leq\exp( - \frac{c}{2} \log{p\choose k} )$.
For\vspace*{1pt} any $\theta\in\calM$, denote by $U\times V$ the support of
$\theta$. Then $\llvert U\rrvert,\llvert V\rrvert
\geq k$.
Let $I,J$ be independently and uniformly drawn at random from all
subsets of cardinality $k$ of $U$ and $V$, respectively. Then $\Expect
[\sum_{i\in I,j\in J} \theta_{ij}] \geq\sum_{i\in U,j\in V} \theta
_{ij}\Expect[{\mathbf{1}_{ \{{i\in I} \}}}{\mathbf
{1}_{ \{{j\in J} \}}}] = \frac{k^2}{\llvert U\rrvert\llvert V\rrvert
} \sum_{i\in
U,j\in V} \theta_{ij} \geq\lambda k^2$.\break Therefore
there exist $S\subset U$ and $T\subset V$, such that $\llvert
S\rrvert=\llvert T\rrvert=k$ and\break 
$\sum_{i \in S, j\in T} \theta_{ij} \geq\lambda k^2$. Then $\sum_{i
\in S, j\in T} X_{ij} \sim N(\mu, k^2)$, where $\mu\geq\lambda k^2$.
Therefore $\mathbb{P}_{\theta} \{ \Tscan\leq\tau' \}
\leq\mathbb{P}_{\theta} \{ \sum_{i \in S, j\in T} X_{ij} \leq
k \tau' \} \leq\exp( -\frac{(\mu-k\tau')_+^2}{2 k}
) \leq\break \exp( -\frac{1}{2} (\lambda k - \tau')_+^2
)$.\vspace*{1pt}

The desired bound (\ref{eqtmax}) on $\Tmax$ follows from analogous
arguments since $\Tmax$ coincides with $\Tscan$ with parameter $k=1$.

\subsection{Proof of Lemma~\texorpdfstring{\protect\ref{lmmnormal-trunc}}{2}}
\label{secproof-trunc}

For the first claim, the marginal density function of $W$ is $f_1$ in
(\ref{eqcond-gauss}). So by definition,
\begin{eqnarray*}
\TV\bigl(\calL(W), N(\mu,1)\bigr) &=& \frac{1}{2} \int
_{\reals} \bigl\llvert f_1(x) - \varphi(x-\mu)\bigr
\rrvert\,\diff x = \overline{\Phi}(M-\mu) +\overline\Phi(M+\mu)
\\
& \leq& 2\overline{\Phi}(M-\mu) \leq\exp\bigl( -{(M-\mu)^2}/{2} \bigr)
\leq\exp\bigl( -\bigl(M^2-1\bigr)/2 \bigr), 
\end{eqnarray*}
were the last inequality is due to the fact that for any $0<\mu\leq
\frac{1}{2M}$, $(M-\mu)^2 \geq(M-\frac{1}{2M})^2 \geq M^2-1$.
For\vspace*{1pt} the second claim, the marginal density function of $W$ is
$f = \frac{1}{2}(f_0 + f_1) = c_0\varphi(x){\mathbf{1}_{ \{
{\llvert x\rrvert\leq M} \}}}$.
Thus
\begin{eqnarray*}
\TV\bigl(\calL(W), N(0,1)\bigr) &=& \frac{1}{2} \int_{\reals}
\bigl\llvert c_0\varphi(x){\mathbf{1}_{ \{{\llvert x\rrvert\leq
M} \}}} - \varphi(x)
\bigr\rrvert\,\diff x = 2\overline\Phi(M)
\\
&\leq&\exp\bigl( -{M^2}/{2} \bigr). 
\end{eqnarray*}
This completes the proof.

\subsection{Proof of Lemma~\texorpdfstring{\protect\ref{lmmnull-reduce}}{3}}
\label{secproof-null-reduce}

We need the following result on the total variation between product
distributions.

%
\begin{lemma}
$\TV( \prod_{i=1}^n P_i, \prod_{i=1}^n Q_i )
\leq\sum_{i=1}^n \TV(P_i,Q_i)$.
\label{lmmtv-prod}
\end{lemma}

\begin{pf}
Recall the dual representation of the total variation \cite{strassenmarginal},
%
\begin{equation}
\TV(P,Q) = \min_{P_{AB}} \bigl\{{ \mathbb{P} \{ A\neq
B \} }\dvtx P_A=P, P_B=Q\bigr\} \label{eqstrassen}
\end{equation}
with infimum over all couplings of $P$ and $Q$. Denote by $P_{A_iB_i}$
the optimal coupling of $P_i$ and $Q_i$ so that ${ \mathbb{P} \{ A_i\neq B_i \} }=\TV(P_i,Q_i)$. Then
$\prod_{i=1}^n P_{A_iB_i}$ is a coupling between the product measures, and
the conclusion follows from the union bound.
\end{pf}

\begin{pf*}{Proof of Lemma~\ref{lmmnull-reduce}}
Let $\widetilde{B} \in\reals^{N_2\times N_2}$ have i.i.d. $N(0,1)$ entries
and be independent of $A$. Let $\widetilde{X} \in\reals^{p\times p}$ be
obtained by applying operations (\ref{eqB-block}) and (\ref
{eqX-reduce}) to $\widetilde{B}$ instead of $B$.
Then it is straightforward to verify that $\widetilde{X}$ has i.i.d. $N(0,1)$
entries, that is, $\calL(\widetilde{X}) = \Prob_0$. Hence
%
\begin{eqnarray}\label{eqBtB-null}
\TV\bigl(\calL(X), \calL(\widetilde{X})\bigr) & \leq&\TV\bigl(\calL(B),
\calL(
\widetilde{B})\bigr)\nonumber
\\
&=& \TV\Biggl( \prod
_{i,j=1}^{N_2} \calL(B_{ij}), \prod
_{i,j=1}^{N_2} N(0,1) \Biggr)
\\
& \leq &\sum_{i,j=1}^{N_2} \TV\bigl(
\calL(B_{ij}), N(0,1)\bigr) 
\leq
N_2^2 \eexp^{-M^2/2} = \frac{1}{4N}.\nonumber
\end{eqnarray}
Here, the first inequality is due to the data processing inequality for
the total variation \cite{Csiszar67}, the second last inequality is
due to Lemma~\ref{lmmtv-prod} and the last inequality is due to
Lemma~\ref{lmmnormal-trunc}.
\end{pf*}

\subsection{Proof of Lemma~\texorpdfstring{\protect\ref{lmmalter-reduce}}{4}}
\label{secproof-alter-reduce}

Recall that $N$ is even with $N_2 = N/2$.
When $A\sim\calG(N,1/2,\kappa)$, let $V \subset[N]$ denote the
vertex subset of size $\kappa$ on which the planted clique in $A$ is supported.
For any subset $S \subset\{N_2+1,\dots, N\}$, we have $S-N_2 \subset[N_2]$.
Further define
%
\begin{equation}
\label{eqV} V_1 = \bigl(V\cap\{N_2+1,\dots, N\}
\bigr)-N_2, \qquad V_2 = V\cap[N_2].
\end{equation}
Then $\llvert V_1\rrvert+\llvert V_2\rrvert =
\kappa$, and the $A_0$ matrix has all ones on
$V_1\times V_2$ and i.i.d. $\mathrm{Bernoulli}(1/2)$ entries elsewhere.
Define $h\dvtx [N_2]\to[p]$ by
%
\begin{equation}
h(x) = 1 + (x-1) \operatorname{mod} p. \label{eqh}
\end{equation}
%
For $i=1,2$, let
%
\begin{equation}
U_i = h(V_i). \label{eqU}
\end{equation}
By the definition of $X$, for each $a,b\in[p]$, we can define sets
%
\begin{eqnarray}
N_{ab}&\triangleq& \bigl[h^{-1}(a)\times
h^{-1}(b)\bigr] \setminus(V_1 \times V_2),
\label{eqNab}
\\
T_{ab}&\triangleq& \bigl[h^{-1}(a)\times
h^{-1}(b)\bigr] \cap(V_1\times V_2).
\label{eqTab}
\end{eqnarray}


$1^\circ$ We first show that the event
%
\begin{equation}
\label{eqE1} E = \bigl\{\llvert U_1\rrvert\geq k\bigr\}\cap\bigl
\{ \llvert U_2\rrvert\geq k\bigr\}
\end{equation}
occurs with high probability.
To this end, first note that
%
\begin{eqnarray}\label{eqV1kappa}
\Prob\bigl\{ \llvert V_1\rrvert< \kappa/4 \bigr\} & \leq&\sum
_{j=1}^{\kappa/4} \frac{{N_2\choose j}{N_2\choose\kappa
-j}}{{N\choose\kappa}} \leq
\frac{\kappa}{4}\frac{{N_2\choose\kappa/4}{N_2\choose3\kappa
/4}}{{N\choose\kappa}} = \frac{\kappa}{4} \frac{{\kappa\choose\kappa
/4}{N-\kappa\choose N_2-\kappa
/4}}{{N\choose N_2}}
\nonumber\\[-8pt]\\[-8pt]\nonumber
& \leq&\frac{\kappa}{4} \biggl( \frac{\eexp}{4} \biggr)^{\kappa/4}
\sqrt{\frac
{2 N}{N-\kappa}} \leq\frac{\kappa}{2\sqrt{2}} \biggl( \frac{\eexp}{4}
\biggr)^{\kappa/4},
\end{eqnarray}
where the second inequality is due to the fact that $j \mapsto
{N_2\choose j}{N_2\choose\kappa-j}$ is increasing for $j \leq(\kappa
-1)/2$, the third inequality is by the bound on the central binomial
coefficient $\frac{2^{2n}}{\sqrt{4n}} \leq\binom{2n}{n} \leq\frac
{2^{2n}}{\sqrt{2n}}$ \cite{Koshy09},\vspace*{-1pt} equation~(2.12), and $\binom
{n}{k} \leq(\frac{\eexp n}{k})^k$ and the last inequality is due to
$N\geq2\kappa$.
By symmetry, since $\llvert V_1\rrvert \stackrel{\mathrm{(d)}}{=}\llvert V_2\rrvert$ and
$\llvert V_1\rrvert+\llvert V_2\rrvert=\kappa$,
$\Prob\{ \llvert V_2\rrvert < \kappa/4 \} =
\Prob\{ \llvert V_1\rrvert >
3\kappa/4 \}$ also
satisfies the upper bound (\ref{eqV1kappa}).

Note that conditioning on the size $\llvert V_1\rrvert =
\kappa_1$, the set $V_1$
is chosen uniformly at random among all $\kappa_1$ subsets of $[N_2]$.
Thus, for any $\kappa_1 \in[\kappa/4, 3\kappa/4]$ and $c_0=1/20$,
%
\begin{eqnarray*}
\Prob\bigl\{ \llvert U_1\rrvert< c_0\kappa\mid
\llvert V_1\rrvert= \kappa_1 \bigr\} &\leq& \sum
_{j= { \lceil{\kappa_1/\ell} \rceil
}}^{c_0\kappa} \frac{ {p\choose
j} {j\ell\choose\kappa_1} }{{N_2\choose\kappa_1}}
\leq  c_0 \kappa\frac{ {p\choose c_0\kappa} {c_0 \kappa\ell
\choose\kappa_1} }{{N_2\choose\kappa_1}}
\\
&\leq& c_0\kappa
\biggl( {\eexp p  \over c_0\kappa}
\biggr)^{c_0\kappa} \biggl( {\eexp c_0 \kappa\ell
\over\kappa_1 }\biggr)^{\kappa_1} \biggl({
\kappa_1  \over N_2 -
\kappa_1 }\biggr)^{\kappa_1}
\\
& \leq& c_0\kappa\exp\biggl( c_0\kappa\log
\frac{\eexp p}{c_0\kappa
} - \frac{\kappa}{4}\log\frac{N_2-\kappa}{\eexp c_0\kappa\ell} \biggr)
\\
&\leq&\frac{\kappa}{20} \exp\biggl( -\frac{\kappa}{5}\log\frac
{p}{\kappa}
\biggr).
\end{eqnarray*}
Here the first inequality is because $j \mapsto{p\choose j} {j\ell
\choose\kappa_1}$ is increasing for $j \leq(p-1)/2$,
the last inequality holds under the assumption that $\kappa\ge20$ and
$p\geq2\kappa$.
Since $k={ \lfloor{\kappa/20} \rfloor}$, the last two
displays together lead to
\begin{eqnarray*}
\Prob\bigl\{ \llvert U_1\rrvert< k \bigr\} &\leq& \sum
_{\kappa_1 = 0}^\kappa\Prob\bigl\{ \llvert U_1
\rrvert< \kappa/20 \mid\llvert V_1\rrvert= \kappa_1
\bigr\} \Prob\bigl\{ \llvert V_1\rrvert= \kappa_1
\bigr\}
\\
&\leq&\Prob\bigl\{ \llvert V_1\rrvert< \kappa/4 \bigr\} + \Prob
\bigl\{ \llvert V_1\rrvert> 3\kappa/4 \bigr\}
\\
&&{} + \max
_{\kappa_1 \in[\kappa/4, 3\kappa/4]}\Prob\bigl\{ \llvert U_1\rrvert<
\kappa/20
\mid\llvert V_1\rrvert= \kappa_1 \bigr\}
\\
&\leq& \kappa\biggl( \frac{\eexp}{4} \biggr)^{\kappa/4} +
\frac
{\kappa}{20} \exp\biggl( -\frac{\kappa}{5}\log\frac{p}{\kappa} \biggr),
\end{eqnarray*}
and the union bound further leads to
%
\begin{equation}
\begin{aligned} \label{eqprob-Ec} \Prob\bigl\{ E^c \bigr\}
& \leq& 2\Prob\bigl\{ \llvert U_1\rrvert< k \bigr\} \leq2 \kappa
\biggl( \frac{\eexp}{4} \biggr)^{\kappa/4} + \frac
{\kappa}{10} \exp
\biggl( -\frac{\kappa}{5}\log\frac{p}{\kappa} \biggr)
\\
& \leq& 40k \biggl( \frac{\eexp}{4} \biggr)^{5k} + 2k \exp\biggl(
-4k\log\frac{p}{20k} \biggr).
\end{aligned}
\end{equation}

$2^\circ$
Conditioned on $V$, we generate a random matrix $\widetilde{B} =
(\widetilde
{B}_{ij})\in\reals^{N_2\times N_2}$ with independent entries where
%
\begin{eqnarray}
&&\qquad \widetilde{B}_{ij} \sim N(\mu,1) \qquad\mbox{if } (i,j)\in
V_1\times V_2, \qquad\widetilde{B}_{ij}\sim
N(0,1)\qquad\mbox{otherwise}. \label{eqtB}
\end{eqnarray}
%
Then we apply (\ref{eqB-block}) and (\ref{eqX-reduce}) to
$\widetilde{B}$
instead of $B$ to obtain a $p\times p$ random matrix~$\widetilde{X}$.
The intuition is that $\tB$ and $\tX$ correspond to the ideal input
and output of the reduction scheme, in the sense that the $\calL(\tX
)$ is, as we show next, close to a desired mixture on the alternative
hypotheses. Our choice of the distribution $\calF_0$ and $\calF_1$
ensures that $B$ is close to the ideal case $\tB$ in total variation,
and the data processing inequality guarantees that the output $X$ is
also close to $\tX$.

To this end, note that conditioned on $V$, for each $a,b\in[p]$, we have
%
\begin{equation}
\widetilde{X}_{ab} = \frac{1}{\ell} \sum
_{i \in h^{-1}(a), j\in h^{-1}(b)} \widetilde{B}_{ij} = \frac{1}{\ell}
\biggl( \sum_{(i,j)\in N_{ab}} \widetilde{B}_{ij} +
\sum_{(i,j)\in T_{ab}} \widetilde{B}_{ij} \biggr),
\label{eqtX}
\end{equation}
where the sets $N_{ab}$ and $T_{ab}$ are defined in (\ref{eqNab}) and
(\ref{eqTab}), respectively.
The last two displays together imply that $\widetilde{X}_{ab}$ follows the
Gaussian distribution with unit variance and mean $\Expect[\widetilde
{X}_{ab}] = \frac{\mu\llvert T_{ab}\rrvert}{\ell}$.
Since for any $(a,b)\in(U_1, U_2)$, $\llvert T_{ab}\rrvert
\geq1$, we have $\Expect
[\widetilde{X}_{ab}]
\geq\frac{\mu}{\ell} = \frac{2\mu p}{N}$.
Last but not the least, since the entries of $\widetilde{X}$ are sums of
mutually independent random variables, they are mutually independent themselves.
Note that for each fixed $V$ [and hence fixed $(V_1,V_2)$ and
$(U_1,U_2)$], $\Indc_{E}$ is deterministic.
Therefore, for any $V$ such that $\Indc_E=1$, there exists some
$\theta= \theta(V) \in\calM(p,k,\frac{2\mu p}{N})$ such that
$\calL(\widetilde{X}\mid V) = \Prob_\theta$.
Define the probability distribution $\pi= \calL(\theta(V)\mid E)$, which
is supported on the set $\calM(p,k,\frac{2\mu p}{N})$.
Then $\calL(\tX\mid E) = \Prob_\pi$, which is a mixture of
distributions of $\{\Prob_\theta\dvtx \theta\in\calM(p,k,\frac{2\mu
p}{N})\}$ with respect to the prior $\pi$.

It remains to show that the law of $X$ is close to the mixture $\Prob
_\pi$.
By the convexity of the $(P,Q)\mapsto\TV(P,Q)$, we have
%
\begin{eqnarray} \label{eqtv-bd-V}
\TV\bigl(\calL(X), \calL(\widetilde{X})\bigr) &\leq& \Expect_V\bigl[
\TV\bigl(\calL(X\mid V), \calL(\widetilde{X}\mid V)\bigr)\bigr]\nonumber
\\
& \leq& \Expect_V\bigl[\TV\bigl(\calL(B\mid V), \calL(\widetilde{B}\mid V)
\bigr)\bigr]
\nonumber\\[-8pt]\\[-8pt]\nonumber
&\leq& \sum_{i,j=1}^{N_2} \TV\bigl(
\calL(B_{ij}\mid V), \calL(\widetilde{B}_{ij}\mid V)\bigr)
\\
&\leq& N_2^2 \eexp^{(M^2-1)/2} \leq\frac{\eexp^{1/2}}{4N}
\leq\frac{1}{p},\nonumber
\end{eqnarray}
where
the second inequality is by the data processing inequality,
the third inequality is due to Lemma~\ref{lmmtv-prod} since
conditioned on $V$ both $B$ and $\tB$ have independent entries, and
the fourth inequality is by Lemma~\ref{lmmnormal-trunc}, and the last
inequality follows from the assumption that $M \geq\sqrt{6\log N}$.
Finally, using $\TV(\calL(\tX),\calL(\tX\mid E)) = \mathbb{P} \{ {E^{\mathrm{c}}} \} $, we obtain
%
\begin{eqnarray}
\TV\bigl(\calL(X), \bbP_\pi\bigr) & \leq&\TV\bigl(\calL(X), \calL(
\tX)\bigr)+\TV\bigl(\calL(\tX), \calL(\tX\mid E)\bigr) \label
{eqcont-tv-bd-decomp}
\nonumber\\[-8pt]\\[-8pt]\nonumber
& \leq&\frac{1}{p} + \Prob\bigl\{ E^c \bigr\},
\end{eqnarray}
where the last inequality is due to (\ref{eqtv-bd-V}).
In view of
(\ref{eqprob-Ec}), this completes the proof.


\begin{supplement}[id=suppA]
\stitle{Supplement to ``Computational barriers in minimax submatrix detection''}
\slink[doi]{10.1214/14-AOS1300SUPP} 
\sdatatype{.pdf}
\sfilename{aos1300\_supp.pdf}
\sdescription{We provide proofs of Theorem \ref{thmequiv} and Lemmas \ref{lmmnull-discrete} and \ref{lmmalter-discrete}.}
\end{supplement}

%

\printaddresses

\begin{thebibliography}{42}
\bibitem{Addario10}
%
\begin{barticle}[mr]
\bauthor{\bsnm{Addario-Berry},~\bfnm{Louigi}\binits{L.}},
\bauthor{\bsnm{Broutin},~\bfnm{Nicolas}\binits{N.}},
\bauthor{\bsnm{Devroye},~\bfnm{Luc}\binits{L.}} \AND
\bauthor{\bsnm{Lugosi},~\bfnm{G{\'a}bor}\binits{G.}}
(\byear{2010}).
\btitle{On combinatorial testing problems}.
\bjournal{Ann. Statist.}
\bvolume{38}
\bpages{3063--3092}.
\bid{doi={10.1214/10-AOS817}, issn={0090-5364}, mr={2722464}}
\end{barticle}
%

\bptok{imsref}%
\endbibitem

\bibitem{AAK07}
%
\begin{bincollection}[mr]
\bauthor{\bsnm{Alon},~\bfnm{Noga}\binits{N.}},
\bauthor{\bsnm{Andoni},~\bfnm{Alexandr}\binits{A.}},
\bauthor{\bsnm{Kaufman},~\bfnm{Tali}\binits{T.}},
\bauthor{\bsnm{Matulef},~\bfnm{Kevin}\binits{K.}},
\bauthor{\bsnm{Rubinfeld},~\bfnm{Ronitt}\binits{R.}} \AND
\bauthor{\bsnm{Xie},~\bfnm{Ning}\binits{N.}}
(\byear{2007}).
\btitle{Testing {$k$}-wise and almost {$k$}-wise independence}.
In \bbooktitle{S{TOC}'07---{P}roceedings of the 39th {A}nnual ACM
{S}ymposium on {T}heory of {C}omputing}
\bpages{496--505}.
\bpublisher{ACM},
\blocation{New York}.
\bid{doi={10.1145/1250790.1250863}, mr={2402475}}
\end{bincollection}
%

\bptok{imsref}%
\endbibitem

\bibitem{Alon98}
%
\begin{binproceedings}[mr]
\bauthor{\bsnm{Alon},~\bfnm{Noga}\binits{N.}},
\bauthor{\bsnm{Krivelevich},~\bfnm{Michael}\binits{M.}} \AND
\bauthor{\bsnm{Sudakov},~\bfnm{Benny}\binits{B.}}
(\byear{1998}).
\btitle{Finding a large hidden clique in a random graph}.
In \bbooktitle{Proceedings of the {N}inth {A}nnual ACM--SIAM
{S}ymposium on {D}iscrete {A}lgorithms ({S}an {F}rancisco, CA, 1998)}
\bpages{594--598}.
\bpublisher{ACM},
\blocation{New York}.
\bid{mr={1642973}}
\end{binproceedings}
%

\bptok{imsref}%
\endbibitem

\bibitem{Ames11}
%
\begin{barticle}[mr]
\bauthor{\bsnm{Ames},~\bfnm{Brendan~P.~W.}\binits{B.~P.~W.}} \AND
\bauthor{\bsnm{Vavasis},~\bfnm{Stephen~A.}\binits{S.~A.}}
(\byear{2011}).
\btitle{Nuclear norm minimization for the planted clique and biclique
problems}.
\bjournal{Math. Program.}
\bvolume{129}
\bpages{69--89}.
\bid{doi={10.1007/s10107-011-0459-x}, issn={0025-5610}, mr={2831403}}
\end{barticle}
%

\bptok{imsref}%
\endbibitem

\bibitem{ABW10}
%
\begin{bincollection}[mr]
\bauthor{\bsnm{Applebaum},~\bfnm{Benny}\binits{B.}},
\bauthor{\bsnm{Barak},~\bfnm{Boaz}\binits{B.}} \AND
\bauthor{\bsnm{Wigderson},~\bfnm{Avi}\binits{A.}}
(\byear{2010}).
\btitle{Public-key cryptography from different assumptions}.
In \bbooktitle{S{TOC}'10---{P}roceedings of the 2010 ACM
{I}nternational {S}ymposium on {T}heory of {C}omputing}
\bpages{171--180}.
\bpublisher{ACM},
\blocation{New York}.
\bid{mr={2743266}}
\end{bincollection}
%

\bptok{imsref}%
\endbibitem

\bibitem{Arias13a}
%
\begin{bmisc}[auto:parserefs-M02]
\bauthor{\bsnm{Arias-Castro},~\bfnm{E.}\binits{E.}} \AND
\bauthor{\bsnm{Verzelen},~\bfnm{N.}\binits{N.}}
(\byear{2013}).
\bhowpublished{Community detection in random networks.
Preprint. Available at \arxivurl{arXiv:1302.7099}.}
\end{bmisc}
%

\bptok{imsref}%
\endbibitem

\bibitem{Arora09}
%
\begin{bbook}[mr]
\bauthor{\bsnm{Arora},~\bfnm{Sanjeev}\binits{S.}} \AND
\bauthor{\bsnm{Barak},~\bfnm{Boaz}\binits{B.}}
(\byear{2009}).
\btitle{Computational Complexity: A Modern Approach}.
\bpublisher{Cambridge Univ. Press},
\blocation{Cambridge}.
\bid{doi={10.1017/CBO9780511804090}, mr={2500087}}
\end{bbook}
%

\bptok{imsref}%
\endbibitem

\bibitem{BKRSW11}
%
\begin{bincollection}[auto:parserefs-M02]
\bauthor{\bsnm{Balakrishnan},~\bfnm{S.}\binits{S.}},
\bauthor{\bsnm{Kolar},~\bfnm{M.}\binits{M.}},
\bauthor{\bsnm{Rinaldo},~\bfnm{A.}\binits{A.}},
\bauthor{\bsnm{Singh},~\bfnm{A.}\binits{A.}} \AND
\bauthor{\bsnm{Wasserman},~\bfnm{L.}\binits{L.}}
(\byear{2011}).
\btitle{Statistical and computational tradeoffs in biclustering}.
In \bbooktitle{NIPS 2011 Workshop on Computational Trade-Offs in
Statistical Learning}.
\end{bincollection}
%

\bptok{imsref}%
\endbibitem

\bibitem{Berthet13}
%
\begin{barticle}[auto:parserefs-M02]
\bauthor{\bsnm{Berthet},~\bfnm{Q.}\binits{Q.}} \AND
\bauthor{\bsnm{Rigollet},~\bfnm{P.}\binits{P.}}
(\byear{2013}).
\btitle{Complexity theoretic lower bounds for sparse principal
component detection}.
\bjournal{Journal of Machine Learning Research: Workshop and
Conference Proceedings}
\bvolume{30}
\bpages{1--21}.
\end{barticle}
%

\bptok{imsref}%
\endbibitem

\bibitem{Berthet12}
%
\begin{barticle}[mr]
\bauthor{\bsnm{Berthet},~\bfnm{Quentin}\binits{Q.}} \AND
\bauthor{\bsnm{Rigollet},~\bfnm{Philippe}\binits{P.}}
(\byear{2013}).
\btitle{Optimal detection of sparse principal components in high dimension}.
\bjournal{Ann. Statist.}
\bvolume{41}
\bpages{1780--1815}.
\bid{doi={10.1214/13-AOS1127}, issn={0090-5364}, mr={3127849}}
\end{barticle}
%

\bptok{imsref}%
\endbibitem

\bibitem{Bhamidi12}
%
\begin{bmisc}[auto:parserefs-M02]
\bauthor{\bsnm{Bhamidi},~\bfnm{S.}\binits{S.}},
\bauthor{\bsnm{Dey},~\bfnm{P.~S.}\binits{P.~S.}} \AND
\bauthor{\bsnm{Nobel},~\bfnm{A.~B.}\binits{A.~B.}}
(\byear{2012}).
\bhowpublished{Energy landscape for large average submatrix detection
problems in {G}aussian random matrices.
Preprint. Available at \arxivurl{arXiv:1211.2284}.}
\end{bmisc}
%

\bptok{imsref}%
\endbibitem

\bibitem{Bhatia}
%
\begin{bbook}[mr]
\bauthor{\bsnm{Bhatia},~\bfnm{Rajendra}\binits{R.}}
(\byear{1997}).
\btitle{Matrix Analysis}.
\bseries{Graduate Texts in Mathematics}
\bvolume{169}.
\bpublisher{Springer},
\blocation{New York}.
\bid{doi={10.1007/978-1-4612-0653-8}, mr={1477662}}
\end{bbook}
%

\bptok{imsref}%
\endbibitem

\bibitem{BI12}
%
\begin{barticle}[mr]
\bauthor{\bsnm{Butucea},~\bfnm{Cristina}\binits{C.}} \AND
\bauthor{\bsnm{Ingster},~\bfnm{Yuri~I.}\binits{Y.~I.}}
(\byear{2013}).
\btitle{Detection of a sparse submatrix of a high-dimensional noisy matrix}.
\bjournal{Bernoulli}
\bvolume{19}
\bpages{2652--2688}.
\bid{doi={10.3150/12-BEJ470}, issn={1350-7265}, mr={3160567}}
\end{barticle}
%

\bptok{imsref}%
\endbibitem

\bibitem{Candes09}
%
\begin{barticle}[mr]
\bauthor{\bsnm{Cand{\`e}s},~\bfnm{Emmanuel~J.}\binits{E.~J.}} \AND
\bauthor{\bsnm{Recht},~\bfnm{Benjamin}\binits{B.}}
(\byear{2009}).
\btitle{Exact matrix completion via convex optimization}.
\bjournal{Found. Comput. Math.}
\bvolume{9}
\bpages{717--772}.
\bid{doi={10.1007/s10208-009-9045-5}, issn={1615-3375}, mr={2565240}}
\end{barticle}
%

\bptok{imsref}%
\endbibitem

\bibitem{ChJo13}
%
\begin{barticle}[mr]
\bauthor{\bsnm{Chandrasekaran},~\bfnm{Venkat}\binits{V.}} \AND
\bauthor{\bsnm{Jordan},~\bfnm{Michael~I.}\binits{M.~I.}}
(\byear{2013}).
\btitle{Computational and statistical tradeoffs via convex relaxation}.
\bjournal{Proc. Natl. Acad. Sci. USA}
\bvolume{110}
\bpages{E1181--E1190}.
\bid{doi={10.1073/pnas.1302293110}, issn={1091-6490}, mr={3047651}}
\end{barticle}
%

\bptok{imsref}%
\endbibitem

\bibitem{cover}
%
\begin{bbook}[mr]
\bauthor{\bsnm{Cover},~\bfnm{Thomas~M.}\binits{T.~M.}} \AND
\bauthor{\bsnm{Thomas},~\bfnm{Joy~A.}\binits{J.~A.}}
(\byear{2006}).
\btitle{Elements of Information Theory},
\bedition{2nd} ed.
\bpublisher{Wiley},
\blocation{Hoboken, NJ}.
\bid{mr={2239987}}
\end{bbook}
%

\bptok{imsref}%
\endbibitem

\bibitem{Csiszar67}
%
\begin{barticle}[mr]
\bauthor{\bsnm{Csisz{\'a}r},~\bfnm{I.}\binits{I.}}
(\byear{1967}).
\btitle{Information-type measures of difference of probability
distributions and indirect observations}.
\bjournal{Studia Sci. Math. Hungar.}
\bvolume{2}
\bpages{299--318}.
\bid{issn={0081-6906}, mr={0219345}}
\end{barticle}
%

\bptok{imsref}%
\endbibitem

\bibitem{Dekel10}
%
\begin{bincollection}[mr]
\bauthor{\bsnm{Dekel},~\bfnm{Yael}\binits{Y.}},
\bauthor{\bsnm{Gurel-Gurevich},~\bfnm{Ori}\binits{O.}} \AND
\bauthor{\bsnm{Peres},~\bfnm{Yuval}\binits{Y.}}
(\byear{2011}).
\btitle{Finding hidden cliques in linear time with high probability}.
In \bbooktitle{A{NALCO}11---{W}orkshop on {A}nalytic {A}lgorithmics
and {C}ombinatorics}
\bpages{67--75}.
\bpublisher{SIAM},
\blocation{Philadelphia, PA}.
\bid{mr={2815485}}
\end{bincollection}
%

\bptok{imsref}%
\endbibitem

\bibitem{Deshpande13}
%
\begin{bmisc}[auto:parserefs-M02]
\bauthor{\bsnm{Deshpande},~\bfnm{Y.}\binits{Y.}} \AND
\bauthor{\bsnm{Montanari},~\bfnm{A.}\binits{A.}}
(\byear{2013}).
\bhowpublished{Finding hidden cliques of size $\sqrt{N/e}$ in nearly
linear time.
Preprint. Available at \arxivurl{arXiv:1304.7047}.}
\end{bmisc}
%

\bptok{imsref}%
\endbibitem

\bibitem{Feige00}
%
\begin{barticle}[mr]
\bauthor{\bsnm{Feige},~\bfnm{Uriel}\binits{U.}} \AND
\bauthor{\bsnm{Krauthgamer},~\bfnm{Robert}\binits{R.}}
(\byear{2000}).
\btitle{Finding and certifying a large hidden clique in a semirandom graph}.
\bjournal{Random Structures Algorithms}
\bvolume{16}
\bpages{195--208}.
\bid
{doi={10.1002/(SICI)1098-2418(200003)16:2<195::AID-RSA5>3.3.CO;2-1},
issn={1042-9832}, mr={1742351}}
\end{barticle}
%

\bptok{imsref}%
\endbibitem

\bibitem{Feige10}
%
\begin{bincollection}[mr]
\bauthor{\bsnm{Feige},~\bfnm{Uriel}\binits{U.}} \AND
\bauthor{\bsnm{Ron},~\bfnm{Dorit}\binits{D.}}
(\byear{2010}).
\btitle{Finding hidden cliques in linear time}.
In \bbooktitle{21st {I}nternational {M}eeting on {P}robabilistic,
{C}ombinatorial, and {A}symptotic {M}ethods in the {A}nalysis of
{A}lgorithms ({A}of{A}'10)}
\bpages{189--203}.
\bpublisher{Assoc. Discrete Math. Theor. Comput. Sci.},
\blocation{Nancy}.
\bid{mr={2735341}}
\end{bincollection}
%

\bptok{imsref}%
\endbibitem

\bibitem{Feldman13}
%
\begin{bincollection}[mr]
\bauthor{\bsnm{Feldman},~\bfnm{Vitaly}\binits{V.}},
\bauthor{\bsnm{Grigorescu},~\bfnm{Elena}\binits{E.}},
\bauthor{\bsnm{Reyzin},~\bfnm{Lev}\binits{L.}},
\bauthor{\bsnm{Vempala},~\bfnm{Santosh~S.}\binits{S.~S.}} \AND
\bauthor{\bsnm{Xiao},~\bfnm{Ying}\binits{Y.}}
(\byear{2013}).
\btitle{Statistical algorithms and a lower bound for detecting planted
cliques}.
In \bbooktitle{S{TOC}'13---{P}roceedings of the 2013 ACM {S}ymposium
on {T}heory of {C}omputing}
\bpages{655--664}.
\bpublisher{ACM},
\blocation{New York}.
\bid{doi={10.1145/2488608.2488692}, mr={3210827}}
\end{bincollection}
%

\bptok{imsref}%
\endbibitem

\bibitem{HK11}
%
\begin{barticle}[mr]
\bauthor{\bsnm{Hazan},~\bfnm{Elad}\binits{E.}} \AND
\bauthor{\bsnm{Krauthgamer},~\bfnm{Robert}\binits{R.}}
(\byear{2011}).
\btitle{How hard is it to approximate the best {N}ash equilibrium?}
\bjournal{SIAM J. Comput.}
\bvolume{40}
\bpages{79--91}.
\bid{doi={10.1137/090766991}, issn={0097-5397}, mr={2765712}}
\end{barticle}
%

\bptok{imsref}%
\endbibitem

\bibitem{Jerrum92}
%
\begin{barticle}[mr]
\bauthor{\bsnm{Jerrum},~\bfnm{Mark}\binits{M.}}
(\byear{1992}).
\btitle{Large cliques elude the {M}etropolis process}.
\bjournal{Random Structures Algorithms}
\bvolume{3}
\bpages{347--359}.
\bid{doi={10.1002/rsa.3240030402}, issn={1042-9832}, mr={1179827}}
\end{barticle}
%

\bptok{imsref}%
\endbibitem

\bibitem{JP00}
%
\begin{barticle}[mr]
\bauthor{\bsnm{Juels},~\bfnm{Ari}\binits{A.}} \AND
\bauthor{\bsnm{Peinado},~\bfnm{Marcus}\binits{M.}}
(\byear{2000}).
\btitle{Hiding cliques for cryptographic security}.
\bjournal{Des. Codes Cryptogr.}
\bvolume{20}
\bpages{269--280}.
\bid{doi={10.1023/A:1008374125234}, issn={0925-1022}, mr={1779310}}
\end{barticle}
%

\bptok{imsref}%
\endbibitem

\bibitem{taocp2}
%
\begin{bbook}[mr]
\bauthor{\bsnm{Knuth},~\bfnm{Donald~E.}\binits{D.~E.}}
(\byear{1969}).
\btitle{The Art of Computer Programming. {V}ol. 2: {S}eminumerical Algorithms}.
\bpublisher{Addison-Wesley},
\blocation{Reading, MA}.
\bid{mr={0286318}}
\bptnote{check year}%
\end{bbook}
%

\bptok{imsref}%
\endbibitem

\bibitem{KZ12}
%
\begin{barticle}[mr]
\bauthor{\bsnm{Koiran},~\bfnm{Pascal}\binits{P.}} \AND
\bauthor{\bsnm{Zouzias},~\bfnm{Anastasios}\binits{A.}}
(\byear{2014}).
\btitle{Hidden cliques and the certification of the restricted
isometry property}.
\bjournal{IEEE Trans. Inform. Theory}
\bvolume{60}
\bpages{4999--5006}.
\bid{doi={10.1109/TIT.2014.2331341}, issn={0018-9448}, mr={3245368}}
\bptnote{check year}%
\end{barticle}
%

\bptok{imsref}%
\endbibitem

\bibitem{KBRS11}
%
\begin{barticle}[auto:parserefs-M02]
\bauthor{\bsnm{Kolar},~\bfnm{M.}\binits{M.}},
\bauthor{\bsnm{Balakrishnan},~\bfnm{S.}\binits{S.}},
\bauthor{\bsnm{Rinaldo},~\bfnm{A.}\binits{A.}} \AND
\bauthor{\bsnm{Singh},~\bfnm{A.}\binits{A.}}
(\byear{2011}).
\btitle{Minimax localization of structural information in large noisy
matrices}.
\bjournal{Adv. Neural Inf. Process. Syst.}
\bvolume{24}
\bpages{909--917}.
\end{barticle}
%

\bptok{imsref}%
\endbibitem

\bibitem{Koshy09}
%
\begin{bbook}[mr]
\bauthor{\bsnm{Koshy},~\bfnm{Thomas}\binits{T.}}
(\byear{2009}).
\btitle{Catalan Numbers with Applications}.
\bpublisher{Oxford Univ. Press},
\blocation{Oxford}.
\bid{mr={2526440}}
\end{bbook}
%

\bptok{imsref}%
\endbibitem

\bibitem{KNV13}
%
\begin{bmisc}[auto:parserefs-M02]
\bauthor{\bsnm{Krauthgamer},~\bfnm{R.}\binits{R.}},
\bauthor{\bsnm{Nadler},~\bfnm{B.}\binits{B.}} \AND
\bauthor{\bsnm{Vilenchik},~\bfnm{D.}\binits{D.}}
(\byear{2013}).
\bhowpublished{Do semidefinite relaxations really solve sparse PCA?
Preprint. Available at \arxivurl{arXiv:1306.3690}.}
\end{bmisc}
%

\bptok{imsref}%
\endbibitem

\bibitem{Kuvcera92}
%
\begin{bincollection}[mr]
\bauthor{\bsnm{Ku{{v}{c}}era},~\bfnm{Lud{{v}{e}}k}\binits{L.}}
(\byear{1992}).
\btitle{A generalized encryption scheme based on random graphs}.
In \bbooktitle{Graph-Theoretic Concepts in Computer Science
({F}ischbachau, 1991)}.
\bseries{Lecture Notes in Computer Science}
\bvolume{570}
\bpages{180--186}.
\bpublisher{Springer},
\blocation{Berlin}.
\bid{doi={10.1007/3-540-55121-2_17}, mr={1245056}}
\end{bincollection}
%

\bptok{imsref}%
\endbibitem

\bibitem{Kuvcera95}
%
\begin{barticle}[mr]
\bauthor{\bsnm{Ku{{v}{c}}era},~\bfnm{Lud{{v}{e}}k}\binits{L.}}
(\byear{1995}).
\btitle{Expected complexity of graph partitioning problems}.
\bjournal{Discrete Appl. Math.}
\bvolume{57}
\bpages{193--212}.
\bid{doi={10.1016/0166-218X(94)00103-K}, issn={0166-218X}, mr={1327775}}
\end{barticle}
%

\bptok{imsref}%
\endbibitem

\bibitem{LeCam86}
%
\begin{bbook}[mr]
\bauthor{\bsnm{Le Cam},~\bfnm{Lucien}\binits{L.}}
(\byear{1986}).
\btitle{Asymptotic Methods in Statistical Decision Theory}.
\bpublisher{Springer},
\blocation{New York}.
\bid{doi={10.1007/978-1-4612-4946-7}, mr={0856411}}
\end{bbook}
%

\bptok{imsref}%
\endbibitem

\bibitem{MW13}
%
\begin{bmisc}[auto:parserefs-M02]
\bauthor{\bsnm{Ma},~\bfnm{Z.}\binits{Z.}} \AND
\bauthor{\bsnm{Wu},~\bfnm{Y.}\binits{Y.}}
(\byear{2013}).
\bhowpublished{Volume ratio, sparsity, and minimaxity under unitarily
invariant norms.
Preprint. Available at \arxivurl{arXiv:1306.3609}.}
\end{bmisc}
%

\bptok{imsref}%
\endbibitem

\bibitem{MW13-supp}
%
\begin{bmisc}[author]
\bauthor{\bsnm{Ma},~\bfnm{Z.}\binits{Z.}} \AND
\bauthor{\bsnm{Wu},~\bfnm{Y.}\binits{Y.}}
(\byear{2015}).
\bhowpublished{Supplement to ``Computational barriers in minimax
submatrix detection.''
DOI:\doiurl{10.1214/14-AOS1300SUPP}}.
\bptok{imsref}%
\end{bmisc}
%
\endbibitem
%

\bptok{imsref}%
\endbibitem

\bibitem{renyi}
%
\begin{barticle}[mr]
\bauthor{\bsnm{R{\'e}nyi},~\bfnm{A.}\binits{A.}}
(\byear{1959}).
\btitle{On the dimension and entropy of probability distributions}.
\bjournal{Acta Math. Acad. Sci. Hungar.}
\bvolume{10}
\bpages{193--215 (unbound insert)}.
\bid{issn={0001-5954}, mr={0107575}}
\end{barticle}
%

\bptok{imsref}%
\endbibitem

\bibitem{Rossman10}
%
\begin{bmisc}[mr]
\bauthor{\bsnm{Rossman},~\bfnm{Benjamin}\binits{B.}}
(\byear{2010}).
\bhowpublished{Average-case complexity of detecting cliques.
Ph.D. thesis, Massachusetts Institute of Technology.}
\bid{mr={2873600}}
\end{bmisc}
%

\bptok{imsref}%
\endbibitem

\bibitem{Shabalin09}
%
\begin{barticle}[mr]
\bauthor{\bsnm{Shabalin},~\bfnm{Andrey~A.}\binits{A.~A.}},
\bauthor{\bsnm{Weigman},~\bfnm{Victor~J.}\binits{V.~J.}},
\bauthor{\bsnm{Perou},~\bfnm{Charles~M.}\binits{C.~M.}} \AND
\bauthor{\bsnm{Nobel},~\bfnm{Andrew~B.}\binits{A.~B.}}
(\byear{2009}).
\btitle{Finding large average submatrices in high dimensional data}.
\bjournal{Ann. Appl. Stat.}
\bvolume{3}
\bpages{985--1012}.
\bid{doi={10.1214/09-AOAS239}, issn={1932-6157}, mr={2750383}}
\end{barticle}
%

\bptok{imsref}%
\endbibitem

\bibitem{SS00}
%
\begin{bbook}[mr]
\bauthor{\bsnm{Shiryaev},~\bfnm{A.~N.}\binits{A.~N.}} \AND
\bauthor{\bsnm{Spokoiny},~\bfnm{V.~G.}\binits{V.~G.}}
(\byear{2000}).
\btitle{Statistical Experiments and Decisions: Asymptotic Theory}.
\bseries{Advanced Series on Statistical Science \& Applied Probability}
\bvolume{8}.
\bpublisher{World Scientific},
\blocation{River Edge, NJ}.
\bid{doi={10.1142/9789812779243}, mr={1791434}}
\end{bbook}
%

\bptok{imsref}%
\endbibitem

\bibitem{strassenmarginal}
%
\begin{barticle}[mr]
\bauthor{\bsnm{Strassen},~\bfnm{V.}\binits{V.}}
(\byear{1965}).
\btitle{The existence of probability measures with given marginals}.
\bjournal{Ann. Math. Statist.}
\bvolume{36}
\bpages{423--439}.
\bid{issn={0003-4851}, mr={0177430}}
\end{barticle}
%

\bptok{imsref}%
\endbibitem

\bibitem{Sun13}
%
\begin{barticle}[mr]
\bauthor{\bsnm{Sun},~\bfnm{Xing}\binits{X.}} \AND
\bauthor{\bsnm{Nobel},~\bfnm{Andrew~B.}\binits{A.~B.}}
(\byear{2013}).
\btitle{On the maximal size of large-average and ANOVA-fit submatrices
in a {G}aussian random matrix}.
\bjournal{Bernoulli}
\bvolume{19}
\bpages{275--294}.
\bid{doi={10.3150/11-BEJ394}, issn={1350-7265}, mr={3019495}}
\end{barticle}
%

\bptok{imsref}%
\endbibitem

\bibitem{Verzelen13}
%
\begin{bmisc}[auto:parserefs-M02]
\bauthor{\bsnm{Verzelen},~\bfnm{N.}\binits{N.}} \AND
\bauthor{\bsnm{Arias-Castro},~\bfnm{E.}\binits{E.}}
(\byear{2013}).
\bhowpublished{Community detection in sparse random networks.
Preprint. Available at \arxivurl{arXiv:1308.2955}.}
\end{bmisc}
%

\bptok{imsref}%
\endbibitem
\end{thebibliography}
\end{document}